\documentclass [11pt,twoside]   {article}
\usepackage{amsmath}
\usepackage{amsfonts}
\usepackage{amssymb}

\usepackage{hyperref}

 \newtheorem{lemma}{Lemma}[section]
  \newtheorem{theorem}[lemma]{Theorem}
\newtheorem {example}[lemma]{Example}

  \newtheorem{conjecture}[lemma]{Conjecture}
  \newtheorem{remark}[lemma]{Remark}

\numberwithin{equation}{section}
  \newcommand {\pf}  {\mbox{\sc Proof. \,\,}}
  \newcommand {\qed} {\null \hfill \rule{2mm}{2mm}}

\begin {document}
\title{{\Large{\bf Injectivity of the specialization homomorphism of elliptic curves}}}

\author
{  {Ivica Gusi\'c and Petra Tadi\'c}  \vspace{1ex}\\
}

\date{21.8.2012.}
\maketitle

\begin{abstract}
\noindent  Let $E:y^2=(x-e_1)(x-e_2)(x-e_3)$  be a nonconstant
elliptic curve over $\mathbb{Q}(t)$, where $e_j\in \mathbb{Z}[t]$. We describe a method for
finding a specialization $t\mapsto t_0\in\mathbb{Q}$ such that the
specialization homomorphism is injective. The method can be directly extended  to elliptic curves  with
$e_j\in \mathcal{R}_K[t]$ where $K$ is a number field and $\mathcal{R}_K$ is some UFD such that $\mathcal O_K\subset\mathcal R_K\subset K$. Further, we make a simplification of the method for a special case of quadratic twists. The method is applied to obtain exactly the rank and prove that a set of points are free generators of several elliptic curves over $\mathbb Q(t)$ coming from \cite{Me}.
\end{abstract}


\section{Introduction}
 \label{section1}
  Let
\begin{equation}\label{jedn}
E=E(t):y^2=(x-e_1)(x-e_2)(x-e_3),\ e_j\in \mathbb{Z}[t] .
\end{equation}
be a nonconstant (non-isotrivial) elliptic curve, i.e.  $E$ is not
isomorphic over $\mathbb{Q}(t)$
to an elliptic curve over $\mathbb Q$.
Let $t_0\in\mathbb{Q}$ be such that
$$(e_1-e_2)(e_2-e_3)(e_3-e_1)(t_0)\neq 0.$$ Then the specialization
$E(t_0)$ of $E(t)$ is an elliptic curve over $\mathbb{Q}$. Let
 $\sigma=\sigma_{t_0}: E(\mathbb{Q}(t))\rightarrow
E(t_0)(\mathbb{Q})$ be the corresponding specialization
homomorphism (note that it is well defined). The specialization homomorphism can be defined for general
non-split elliptic surfaces in a more general situation. By the Silverman specialization theorem, it
 is injective for all but finitely many
rational $t_0$. As far as
we know, there is no algorithm for determining such a $t_0$ (for general non-split elliptic surfaces). In this
paper we  improve and extend the method from \cite{G-T}, for
finding a specialization $t\mapsto t_0\in\mathbb{Q}$ such that the
specialization homomorphism is injective, in the case of elliptic
curves  \eqref{jedn}. The improvement leads to an
effective algorithm (see Theorem  \ref{main} and Lemma \ref{effective}).  This algorithm can be directly extended to  number fields $K$ of class number one, where the elliptic curves are of the shape \eqref{jedn} with
$e_j\in \mathcal{O}_K[t]$ (here $\mathcal{O}_K$ is the ring of integers of $K$). In section \ref{section3} we treat the general case, when $K$ is a number field of arbitrary class number.  However, the calculations over general number fields are rather complicated. For example, if the class number of the field $K$ is greater then $1$, the ring of integers $\mathcal{O}_K$ has to be replaced by a suitable UFD (see Theorem \ref{mainh}). In Section \ref{section4} we present a simplification of the extended criterion for quadratic twists $E_g$ of elliptic curves $E:y^2=x^3+ax^2+bx+c,\ a,b,c\in \mathbb{Z}$ with nonconstant polynomials $g$ over $\mathbb{Q}$ (see Theorem \ref{tmtwist}). In Section \ref{section5} we describe and comment a family of quadratic twists coming from Mestre: a family of quadratic twists of the general family of elliptic curves $E=E^{a,b}:y^2=x^3+ax+b$ over $\mathbb{Q}$ with certain $14$th degree polynomials $g=g^{a,b}$ in variable $u$ over $\mathbb{Q}$. It is known that the rank of $E_g$ over $\mathbb{Q}(u)$ is at least $2$ for all $a,b,\ ab\neq 0$. By a general principle, these ranks are at most $6$. In Section \ref{section6} we perform an extensive calculation using our criterion (Theorem \ref{mainh}) for  number fields of class number one (including $\mathbb Q$) and   for a number field of class number two. We prove that the rank is two and that given two points are free generators for a wide class of integers $a,b$. The results suggest that the rank is exactly $2$ and that the certain points $P,Q$ are free generators for all $a,b$.
We used Magma \cite{MAGMA}, Pari \cite{Pari}, and {\tt mwrank} \cite{mwrank} for most of our computations.

This paper has its origins in an idea from the article by professor Andrej Dujella \cite[Theorem 4]{Duj}. We would like to thank him for his kind suggestions and comments.

\section{Elliptic curves $y^2=(x-e_1)(x-e_2)(x-e_3),\ e_j\in \mathbb{Z}[t]$}
\label{section2}
In this section we will work over $\mathbb{Q}$ although all results are valid over arbitrary algebraic number fields $K$ with class number $1$.
Let $E$ be the elliptic curve \eqref{jedn}. We have homorphisms
$\Theta_i:E(\mathbb{Q}(t))\rightarrow \mathbb{Q}(t)^{\times}/(\mathbb{Q}(t)^{\times})^2,\ i=1,2,3$  given by\\
\[\left\{\begin{array}{ll}
\Theta_i(x,y)=(x-e_i)\cdot (\mathbb{Q}(t)^{\times})^2,&\mbox{ if $x\neq e_i$,}\\
\Theta_i(e_i,0)=(e_j-e_i)(e_k-e_i)\cdot (\mathbb{Q}(t)^{\times})^2,&\mbox{ where $i\neq j\neq k\neq i$,}\\
\Theta_i(O)=1\cdot (\mathbb{Q}(t)^{\times})^2,&\mbox{ (here $O$ denotes the neutral element)}.
\end{array}
\right.\]

\begin{lemma}\label{dupla}
  $P\in 2E(\mathbb{Q}(t))$ if and only if $\Theta_i(P)=1\cdot (\mathbb{Q}(t)^{\times})^2$ for $i=1,2,3$.
\end{lemma}
\pf Follows from  \cite{Huse},
Chapter 1, Theorem (4.1), and  Chapter 6,
Proposition (4.3). \qed\medskip

Since $\mathbb{Z}[t]$ is a unique factorization domain (UFD), it is evident that for each $P\in E(\mathbb{Q}(t))$ there exists exactly one
triple $(s_1,s_2,s_3),\  s_i=s_i(P)\in \mathbb{Z}[t], \ i=1,2,3 $, of
non-zero square-free elements from $\mathbb{Z}[t]$, such that
\begin{equation}\label{esi}
\Theta_i(P)= s_i(P)\cdot (\mathbb{Q}(t)^{\times})^2.
\end{equation}
 We will also use notations  $s_i(t)$ for $s_i$. Lemma \ref{dupla} can be reformulated as

\begin{equation}\label{opetdupla}
P\in 2E(\mathbb{Q}(t)),\ {\rm if\ and\ only\ if}\ s_i(P)=1,\ {\rm for}\ i=1,2,3
\end{equation}

 It is
easy to see that
\begin{equation}
 s_1s_2s_3\in \mathbb{Z}[t]^2,
 \end{equation}
 and that, for each $i$ and each prime $p\in \mathbb{Z}[t]$, it must be:
 \begin{equation}\label{ps123}
{\rm if}\ p|s_i\ {\rm then}\ p|s_js_k,\ {\rm where}\
i\neq j\neq k \neq i.
\end{equation}

Let $P\in E(\mathbb{Q}(t))\setminus\{O\}$. Then the first coordinate of $P$ is of the form
\begin{equation}\label{ikspe}
 x(P)=\frac{p(t)}{q(t)^{2}},\  {\rm with}\  p(t),q(t)\in \mathbb{Z}[t]\ {\rm coprime}
 \end{equation}
 (recall that $\mathbb{Z}[t]$ is an UFD). Therefore
\[  \left\{
\begin{array}{c}
p(t)-e_1(t)q^2(t)=s_1(P)\square_{\mathbb Z[t]},\\\
p(t)-e_2(t)q^2(t)=s_2(P)\square_{\mathbb Z[t]},\\
p(t)-e_3(t)q^2(t)=s_3(P)\square_{\mathbb Z[t]},
\end{array} \right.\]
where  $\square_{\mathbb Z[t]}$ denotes a square of an element of
$\mathbb Z[t]$.\\
By this, \eqref{ps123} and the fact that $s_i$ are square-free, we deduce that
\begin{equation}\label{pse}
s_i|(e_j-e_i)(e_k-e_i),\mbox{ where }i\ne j\ne k\ne i
\end{equation}
for each $i$. For example, a prime factor of $s_1$ is also a prime
factor of $s_2s_3$. Assume that it is a prime factor of $s_2$.
Then it is a prime factor of $(e_1-e_2)q^2(t)$, hence it is a
prime factor of $e_1-e_2$.



In the
following theorem we make a refinement of the method from
\cite{G-T}, Theorem 3.2. The proof is a modification of that
proof.

\begin{theorem}\label{main}
Let  $E$ be a nonconstant elliptic curve over $\mathbb Q(t)$, given by the equation
$$E=E(t):y^2=(x-e_1)(x-e_2)(x-e_3), (e_1,e_2,e_3\in \mathbb Z[t]).$$
Let $t_0\in\mathbb Q$ be such that the specialization $E(t_0)$ of $E(t)$ is an elliptic curve. \\
Assume that $t_0$ satisfies the following condition.\\
(A) For every nonconstant square-free divisor $h$ in  $\mathbb Z[t]$ of
$$\mbox{$(e_1-e_2)\cdot (e_1-e_3)$\ \  or \ \ $(e_2-e_1)\cdot(e_2-e_3)$ \ \ or \ \ $(e_3-e_1)\cdot(e_3-e_2)$},$$
the rational number $h(t_0)$ is not a square in $\mathbb{Q}$.\\
Then the specialization homomorphism $\sigma:E(\mathbb{Q}(t))\rightarrow
E(t_0)(\mathbb{Q})$ is injective.
\end{theorem}
\pf
Assume that the conditions of the theorem are satisfied and that $\sigma$ is not an injection. So there exists a point $P\in E(\mathbb{Q}(t))\setminus \{O\}$ such that $\sigma(P)= O$. We will prove that it leads to a contradiction.
First we prove that $P\in 2E(\mathbb{Q}(t)).$  By \eqref{opetdupla}, it is equivalent to proving that  $s_i(t)=1$ for each $i=1,2,3$. Since $\sigma$ is injective on
the torsion part \cite[p. 272--273, proof of Theorem
III.11.4]{Sil}, we may assume that $P\neq (e_i,0),\ i=1,2,3.$ By
$p(t)-e_k(t)q^2(t)=s_k(P)\square_{\mathbb Z[t]}$ and the fact that
$q(t_0)=0$, we get $p(t_0)=s_k(t_0)\square_{\mathbb Q}$. Since
$p(t_0)$ should be a non-zero rational square (recall that $q(t_0)=0$ and $p,q$ are coprime),
we see that $s_i(t_0)$ is a
rational square, for each $i=1,2,3$. We claim that $s_k(t)=1$ for
each $k=1,2,3$, i.e. that $P\in 2E(\mathbb{Q}(t))$.\\
Assume that $s_k(t)$ is non-constant for some $k$.
By the above discussion $s_k(t_0)$ is a rational square, which is in
contradiction with  condition (A) of the theorem (recall that by
\eqref{pse}, $s_k$ is a nonconstant square-free divisor of
$(e_i-e_k)\cdot (e_j-e_k)$ in  $\mathbb Z[t]$, with $i\neq j\neq k\neq i$).
 Therefore  $s_k(t)$ is constant for each $k$. Since $s_k(t)$ is
square-free in ${\mathbb Z[t]}$ and $s_k(t_0)$ is a rational
square, we see that $s_k(t)=1$, for each $k$. It proves that that $P\in 2E(\mathbb{Q}(t))$.\\
We claim that there is $P_1\in E(\mathbb{Q}(t))$ such that $2P_1=P$ and $\sigma (P_1)=O$. Let
$P'_1\in E(\mathbb{Q}(t))$ be any point with $2P'_1=P$. Then
$2\sigma (P'_1)=O$, i.e. $\sigma(P'_1)$ is a $2$-torsion point on the specialized curve. Since $\sigma$ is injective on
the torsion points, there exists a $2$-torsion point $Q\in E(\mathbb{Q}(t))$  such that
$\sigma (Q)=\sigma (P'_1)$. Put $P_1=P'_1-Q$.  Then $2P_1=P$, especially $P_1\neq O$, and $\sigma(P_1)=O$. Note that $P_1$ is of infinite order.  Now
the procedure can be continued with $P_1$ instead of $P$, the contradiction. Therefore $P=O$, i.e. $\sigma$ is injective.
 \qed\medskip

In the following remark we discuss the connection of Theorem \ref{main} and \cite{G-T}, Theorem 3.2.

\begin{remark}
Let $\pm p_1\cdot ...\cdot p_m\cdot f_1\cdot...\cdot f_n,$
be a prime factorization of the square-free part of $(e_1-e_2)\cdot
(e_2-e_3)\cdot (e_3-e_1)$ in  $\mathbb Z[t]$ (here $p_i$ are
rational prime numbers, while $f_j$ are irreducible nonconstant
polynomials from $\mathbb Z[t]$, and we may assume that their leading coefficients are
positive).
 Put $I=\{1,...,m\}$ and
$J=\{1,...,n\}$.
Then the main condition on $t_0$ in \cite{G-T}, Theorem 3.2,  can be
paraphrased as: For each $i\in J$ the integer square-free part of
 $f_i(t_0)$  has at least one prime factor that doesn't appear in the  integer square-free part of any
$f_j(t_0)\ (j\in J,\  j\ne i)$   and doesn't appear in the
factorization  $p_1\cdot ...\cdot p_m$. It is easy to see that
if $t_0$ satisfies this condition, then it satisfies
condition (A) from Theorem \ref{main}, too. The converse is not true. For example, set $e_1=0,\ e_2=t,\ e_3=t^2+10$, hence $(e_1-e_2)\cdot(e_2-e_3)\cdot (e_3-e_1)=t(t^2-t+10)(t^2+10)$. Then $t_0:=2$ satisfies condition (A) from
Theorem \ref{main}, but it does not satisfy condition from \cite{G-T}, Theorem 3.2. Namely, $t(2)=2,\ (t^2-t+10)(2)=12=3\cdot 2^2,\
(t^2+10)(2)=14=2\cdot 7.$
\end{remark}

The following lemma shows that most of integers $t_0$ satisfy
condition (A) from Theorem \ref{main}. It follows from the fact
that curves of genus at least one have finitely many integer
points, and the fact that integer points on genus zero curves are
rare.

\begin{lemma}\label{effective}
Let $\cal T$ denote the set of all integers
$t_0$ that satisfy Condition (A) from Theorem \ref{main}. Then there is an
effectively computable constant $c>0$, such that ${\cal T}\cap
[-c,c]\neq\emptyset$. Therefore, the theorem gives a method for
finding a rational number  $t_0$ such that the specialization
homomorphism $\sigma_{t_0}$ is injective.
\end{lemma}
\pf
Condition (A) in Theorem \ref{main} produces the equations  of the
form $z^2=h(t)$ for certain square-free polynomials $h$ over
$\mathbb{Z}$ of degree $d\geq 1$. If $d\leq 2$, the corresponding
curve has genus $0$, if $d=3$ or $4$ the genus is one, and if
$d\geq 5$ the curve is hyperelliptic with genus $\geq 2$. Recall
that  curves over $\mathbb{Q}$ of genus at least $1$ have only
finitely many integer points. Moreover, for elliptic and hyperelliptic curves,
there are explicit bounds for the height of integer points (\cite{Ba},  \cite{Bu}, Theorem 1; see also
 \cite{E-Si}, Theorem 1 b, for a bound of the number of integer points).  For example, from \cite{Bu}, Theorem 1,
 it follows that for $d\geq 3$ there are effectively computable constants $H=H(h), A=A(h)$ and $c_1=c_1(d)$ such that
 if rational integers $(t,z)$ satisfy $z^2=h(t)$ then
 $$|t|\leq H^2\cdot {\exp}\{c_1\cdot A^{3d^2}\cdot |\Delta_h|^{12 d}\cdot (\log|A\Delta_h|)^{6d^2}\cdot \log\log H\},$$
 where $\Delta_h$ denotes the discriminant of $h$.\\
 If $d=1$ or $d=2$ then
the curve $z^2=h(t)$ may have finitely many or infinitely many integer
points.
 If $d=1$ then there is an effectively computable constant $c_2=c_2(h)$ such that the equation $z^2=h(t)$ has
 $\leq c_2\sqrt{X}$ integer solutions with $|t|\leq X$ for $X\geq 1$. We can take $c_2=2(\sqrt{|a|+|b|}+2)$, for $h(t)=at+b$.\\
 Assume that $d=2$. Let $h(t)=at^2+bt+c$ and let $a=D\cdot k^2$, where $D$ is square-free. Now multiplying $z^2=h(t)$ by $4a$ we get $D(2kz)^2=(2at+b)^2+(4ac-b^2)$. So we see that it is enough to estimate the number of integer solutions for $Dz^2=t^2+B$ where $D$ is a squarefree integer and $B$ a nonzero integer. If $D<0$, then $|t|\leq \sqrt{|B|}$. If $D=1$, then by the unique factorization in $\mathbb{Z}$, then $Dz^2=t^2+B$ has $\leq 2 \tau (B)$ solutions in integers, where $\tau(B)$ denotes the number of positive divisors of $B$.  Finally, if $D\geq 2$ then there is an
 effectively computable constant $c_3=c_3(D,B)$ such that $Dz^2=t^2+B$ has $\leq c_3 \tau (B)\log X$ integer solutions
 with $|t|,|z|\leq X$ for sufficiently large $X$ (see \cite{P-Z}, Lemma 3. for a more
precise estimation).\\
Combining these estimates and the fact that $z^2=h(t)$ has at most two integer solutions with fixed integer value of $t$, we get the statement of the lemma.
   \qed

\section{The case of number fields of arbitrary class number}
\label{section3}

In this section $K$ denotes an algebraic number field with the ring of integers $\mathcal{O}_K$. Here we will generalize the Theorem \ref{main} from $\mathbb Q$ to arbitrary number fields $K$, i.e. to elliptic curves over $K(t)$ given by \eqref{jedn} where $e_j$ are polynomials over a chosen unique factorization domain. For the case $K=\mathbb Q$  the chosen  unique factorization domain  was $\mathbb Z$, for $K$ of class number one it will be $\mathcal O_K$ and for $K$ of  class number at least two it will be a suitable one. 

\begin{remark}\label{h1}
 Let $E$ be
 a nonconstant elliptic curve over $K(t)$ of the shape $$E=E(t):y^2=(x-e_1)(x-e_2)(x-e_3),\ e_j\in \mathcal{O}_K[t].$$
Assume that $K$ has  class number $1$. Then
it is easy to see that the method from Section \ref{section2} and Theorem \ref{main} remains valid if we replace $\mathbb{Q}$ by $K$, $\mathbb Z$ by  $ \mathcal{O}_K$ (which is an UFD here), and $t_0\in \mathbb{Q}$ by $t_0\in K$.
\end{remark}

Generally, (when the class number of $K$ is not necessarily $1$) there exists a unique factorization domain $\mathcal R_K$, $\mathcal O_K\subset \mathcal R_K\subset K$ such that its group of units is finitely generated
(see \cite[p. 94, p. 127]{Knapp} for the description of the construction). For $K$ of class number one we have $\mathcal R_K=\mathcal O_K$, especially for $K=\mathbb Q$ we have $\mathcal R_K=\mathbb Z$.  This fact provides the following
generalization of Theorem \ref{main} and the statement of Remark \ref{h1}.

\begin{theorem}\label{mainh}
Let $K$ be a number field. Let $\mathcal R_K$ be as above a  unique factorization domain such that $\mathcal O_K\subset \mathcal R_K\subset K$ and such that its group of units is finitely generated.  Let  $E$ be a nonconstant elliptic curve over $K(t)$, given by the equation
\begin{equation}\label{geneq}
E=E(t):y^2=(x-e_1)(x-e_2)(x-e_3), (e_1,e_2,e_3\in \mathcal R_K[t]).
\end{equation}
Let $t_0\in K$ be such that the specialization $E(t_0)$ of $E(t)$ is an elliptic curve. \\
Assume that $t_0$ satisfies the following condition.\\
(C) For every nonconstant square-free divisor $h$ in  ${\mathcal {R}}_K[t]$ of
$$\mbox{$(e_1-e_2)\cdot (e_1-e_3)$\ \  or \ \ $(e_2-e_1)\cdot(e_2-e_3)$ \ \ or \ \ $(e_3-e_1)\cdot(e_3-e_2)$},$$
the algebraic number $h(t_0)$ is not a square in $K$.\\
Then the specialization homomorphism $\sigma:E(K(t))\rightarrow
E(t_0)(K)$ is injective.
\end{theorem}
\pf Note that the relations \ref{esi}-\ref{pse} from Section \ref{section2} remain valid after replacing
$\mathbb{Z}[t]$ by $\mathcal R_K[t]$. Now the proof is analogous to the proof of Theorem \ref{main}.
\qed
\medskip

Since the group of invertible elements of $\mathcal R_K$ is finitely generated, to check Condition (C) from Theorem \ref{mainh}
we have to check only finitely many square-free divisors $h$ in  ${\mathcal {R}}_K[t]$.
In Section \ref{section6} we will apply this theorem to a number of examples for  $K$ of class number one (including $K=\mathbb Q$) and for a few  for $K$ of class number two. General implementation of Theorem \ref{mainh} in the case of number fields of class number greater then one requires further investigations. Note also that in this article, we treat only elliptic curves of the shape \eqref{geneq}, so the problem of an extension of our criterion to general elliptic curves
$$y^2=x^3+A(t)x^2+B(t)x+C(t),\ A,B,C\in K(t),$$
remains open (when the equation doesn't factor in the desired form).

It can be seen that there is a variant of Lemma \ref{effective} for elliptic curves \eqref{geneq}.
In the following remark we use another argument to prove that there are a
lot of rational integers $t_0$ satisfying condition (C) from Theorem
\ref{mainh}.
\begin{remark}
According to \cite{Sc}, Section 5, Definition 24,  Theorem 50 and Corollary 1,
for each $F\in\mathbb{C}[z,t]$ either:\\
(i) every congruence class
$\cal C$ in $\mathbb{Z}$ contains a congruence subclass ${\cal C}^{\ast}$ such that
for all $t_0\in{\cal C}^{\ast}$ the polynomial $F(z,t_0)$ has no zero in $K$, or\\
(ii)  $F$ viewed as a polynomial in $z$ has a zero in
$K(t)$.\\ By consecutive applying this to the polynomials
$F[z,t]:=z^2-h(t)$ above, we see that for each congruence class
$\cal C$ in $\mathbb{Z}$ there exists a congruence subclass $\cal C^{\ast}$
of $\cal C$, such that the conditions from Theorem \ref{mainh}  are
satisfied for all $t_0\in \cal C^{\ast}$.
\end{remark}

\section{Nonconstant quadratic twists of elliptic curves $y^2=x^3+ax^2+bx+c,\ a,b,c\in \mathbb{Z}$}
\label{section4}
Let $K$ be a finite extension of $\mathbb{Q}$ with ring of
integers $\mathcal{O}_K$. It is well-known that $\mathcal{O}_K$ is
a UFD (unique factorization domain) if and only if it is a principal
ideal domain, or equivalently, if  the class number of $K$ equals
to $1$. In this section $K$ will always denote the splitting  field
 of a
separable cubic polynomial
$$f(x)=x^3+ax^2+bx+c,\ a,b,c\in \mathbb{Z},$$
especially $K$ is Galois. It is easy to see that either $K=\mathbb{Q}$, $ K$ is a quadratic field over $\mathbb{Q}$,  $K$ is a cubic field over $\mathbb{Q}$ with cyclic Galois group, or $K$ is a sextic field over $\mathbb{Q}$ with the Galois group isomorphic to the symmetric group $\mathbf{S}_3$.
 We will always assume that $K$ has class number $1$. For a domain $\mathcal{A}$ and nonzero elements $u,v\in\mathcal{A}$, we will say that $u,v$ are associate if there exists a unit $\epsilon\in\mathcal{A}$ (i.e. an invertible element) such that $v=\epsilon u$.\\
 Assume first that $K$ is a quadratic number field. Then for a rational prime $H$, either\\
(i)  $H$ remains prime in $\mathcal{O}_K$, or\\
(ii)  $H=\epsilon \mathcal{P}^2$ where $\mathcal{P}$ is a prime in $\mathcal{O}_K$ and $\epsilon$ is a unit in $\mathcal{O}_K$ ($H$ ramifies in $K$), or\\
 (iii)   $H=\pm \mathcal{P}\cdot \bar {\mathcal{P}}$  where $\mathcal{P}$ is a prime in $\mathcal{O}_K$,
 $\bar {\mathcal{P}}$ is the conjugate of $\mathcal{P}$ and $\bar {\mathcal{P}}, \mathcal{P}$ are non-associate ($H$ splits in $K$).\\
Similarly, if $H$ is a non-constant irreducible polynomial from $\mathbb{Z}[x]$, then either\\
(I) $H$ remains irreducible in $\mathcal{O}_K[x]$, or\\
(II) $H=\pm \mathcal{P}\cdot \bar {\mathcal{P}}$  where $\mathcal{P}$ is irreducible in $\mathcal{O}_K[x]$, with
 $\mathcal{P},\bar {\mathcal{P}}$ non-associate ($H$ splits).\\
Assume now that $K$ is a cubic cyclic field, i.e. the discriminant
$D$ of $f$  is a rational square. Let $\tau$ denote a
non-trivial automorphism of $K$.
Then by the decomposition of prime ideals in Galois extensions for a rational prime $H$, either\\
(i)  $H$ remains prime in $\mathcal{O}_K$, or\\
(ii)  $H=\epsilon \mathcal{P}^3$ where $\mathcal{P}$ is a prime in $\mathcal{O}_K$ and $\epsilon$ is a unit in $\mathcal{O}_K$ ($H$ ramifies in $K$), or\\
 (iii)   $H=\pm \mathcal{P}\cdot \mathcal{P}^{\tau}\cdot \mathcal{P}^{\tau^2}$  where $\mathcal{P}$ is a prime in
 $\mathcal{O}_K$, with $\mathcal{P}, \mathcal{P}^{\tau}$ non-associate ($H$ splits in $K$).\\
Similarly, if $H$ is a non-constant irreducible polynomial from $\mathbb{Z}[x]$, then either\\
(I) $H$ remains irreducible in $\mathcal{O}_K[x]$, or\\
(II) $H=\pm \mathcal{P}\cdot \mathcal{P}^{\tau}\cdot \mathcal{P}^{\tau^2}$  where $\mathcal{P}$ is irreducible in $\mathcal{O}_K[x]$ with $\mathcal{P}, \mathcal{P}^{\tau}$ non-associate ($H$ splits).\\
In Theorem \ref{tmtwist} below $E$ is the elliptic curve given by
$$E:y^2=x^3+ax^2+bx+c,\ a,b,c\in \mathbb{Z},$$
$g$ is a nonconstant square free polynomial from $\mathbb{Z}[t]$,
and
$$E_g:y^2=x^3+agx^2+bg^2x+cg^3$$ is the quadratic twist of $E$ with $g$. Recall that we assume that the splitting field $K$ of $f$ has class number $1$.  It is easy to see that in this setting we may apply Theorem \ref{main} directly (see Remark \ref{h1}).  Theorem \ref{tmtwist} enables us to avoid  the calculation in algebraic number fields (when the splitting field is quadratic), or to simplify
it (when the splitting field is cubic with cyclic Galois group). The theorem does not treat the case when the Galois group of $f$ is the symmetric group $\mathbf{S}_3$. First we will prove two lemmas. For a rational prime $p$  we let $v_p$ denote the discrete valuation of $\mathbb{Q}$ at $p$.

\begin{lemma}\label{kvadrat}
Let $d\neq 1$ be a squarefree integer, and let $K=\mathbb{Q}(\sqrt{d})$ be the corresponding quadratic field. Let $D$ denote the discriminant of $K$. Assume that $K$ has class number $1$. Then for an arbitrary nonzero rational number $r$ the following statements are equivalent.
\begin{itemize}
\item (i) There exists a unit $\epsilon\in\mathcal{O}_K$ such that $\epsilon r\in K^2$.

\item (ii) For each rational prime $p$, if $v_p(r)$ is odd then $p|D$.

\item (iii) There exists a divisor $d'$ of $D$ in $\mathbb{Z}$ such that $d'r\in \mathbb{Q}^2.$

\end{itemize}

\end{lemma}
\pf (i) implies (ii). Assume that $v_p(r)$ is odd for a rational prime $p$. Since there is a unit $\epsilon$ in $\mathcal{O}_K$ such that $\epsilon r\in K^2$ we conclude that $p$ ramifies in $K$, which is equivalent with $p|D$.\\
(ii) implies (iii). Let $d''$ be the product of all positive rational primes $p$ such that $v_p(r)$ is odd. Then
$d''r\in \mathbb{Q}^2$ or $-d''r\in \mathbb{Q}^2$, and further $d''|D$.\\
(iii) implies (i) follows from $d'r\in \mathbb{Q}^2,\ d'|D$ and the fact that for each $p|D$ there exist $x_p\in K$ and a unit $\epsilon\in \mathcal{O}_K$ such that $p=\epsilon_px_p^2$. \qed\medskip

\begin{lemma}\label{duplasto}
Let $k$ be a field of characteristic zero with an algebraic closure $\bar k$, and let $E:y^2=x^3+ax^2+bx+c$ be an elliptic curve over $k$. Assume that $x^3+ax^2+bx+c=(x-e_1)(x-e_2)(x-e_3)$ over $\bar k$. Let $P(u,v)\in E(k)$ be such that $u-e_1\in k(e_1)^2,\ u-e_2\in k(e_2)^2,\ u-e_3\in k(e_3)^2$. Then there exists $Q\in E(k)$ such that $2Q=P$.
\end{lemma}
\pf Put $K:=k(e_1,e_2,e_3)$.   By \cite{Huse},
Chapter 1, Theorem (4.1), there exists $Q'\in E(K)$ such that $2Q'=P$. We have to prove that there exists $Q\in E(k)$ such that $2Q=P$.
If $K=k$, it is obvious.\\ Assume that $K$ is a quadratic field over $k$, say $K=k(e_2)=k(e_3)$. Let $\tau$ denote the nontrivial automorphism of $K$ over $k$. Then $e_3=\tau (e_2)$.   If $Q'\notin E(k)$ then $Q'^{\tau}=Q'+(e_1,0)$ (namely, if, for example,  $Q'^{\tau}=Q'+(e_2,0)$, then
$Q'=Q'^{\tau}+(e_3,0)=Q'+(e_1,0)$, a contradiction). Therefore, $(Q'+(e_2,0))^{\tau}=Q'+(e_1,0)+(e_3,0)=Q'+(e_2,0)$, hence
$Q:=Q'+(e_2,0)\in E(k)$ and $2Q=P$, as we need.\\
Assume now, that $K:=k(e_1,e_2,e_3)$ is a cubic cyclic field over $k$. Let $\tau$ denote the  automorphism of $K$ over $k$  such that $\tau(e_1)=e_2$.
If $Q'\notin E(k)$ then we may assume that $Q'^{\tau}=Q'+(e_1,0)$. Hence, $(Q'+(e_2,0))^{\tau}=Q'+(e_1,0)+(e_3,0)=Q'+(e_2,0)$, and we may proceed as above.\\
Assume, finally, that $K:=k(e_1,e_2,e_3)$ is a sextic field over $k$. Let $K_0$ be the quadratic field over $k$ such that $K$ is cubic cyclic over $K_0$. By repeating above argument, we first see that there is $Q\in E(K_0)$ such that $2Q=P$, and after that, that  $Q\in E(k)$.\qed \medskip

\begin{theorem}\label{tmtwist}
Let $f(x):=x^3+ax^2+bx+c,\ a,b,c\in \mathbb{Z},$ be a polynomial without repeated roots, and let $g=g(t)$ be a nonconstant polynomial over $\mathbb{Z}$. Set $E:y^2=x^3+ax^2+bx+c$ and $E_g:y^2=x^3+ag(t)x^2+bg(t)^2x+cg(t)^3$.
Let $t_0\in\mathbb Q$ be such that the specialization $E_g(t_0)$ of $E_g$ is well defined and let $\sigma_{t_0}:E_g(\mathbb Q(t))\rightarrow E_g(t_0)(\mathbb Q)$ be the corresponding specialization.\\
(i) Assume that $c=0$ and that $x^2+ax+b=(x-\theta)(x-\bar\theta)$
is $\mathbb{Q}$-irreducible, with splitting field $K$ having class number $1$. Write $a^2-4b=e^2\cdot d$, where
$d\in \mathbb{Z}$ is square free. Let $D$ denote the discriminant of $K$. Assume
that $t_0$ satisfies conditions\\
(A1) For each nonconstant square free divisor $h$ of $bg$ in
$\mathbb{Z}[t]$, $h(t_0)$ is not a square from $\mathbb{Q}$.\\
(A2) For each nonconstant square free divisor $h$  of $eg$ in
$\mathbb{Z}[t]$, and each square free divisor $d'$ of $D$ in $\mathbb{Z}$,  $d'\cdot h(t_0)$ is not a square from $\mathbb{Q}$.\\
Then the specialization homomorphism $\sigma_{t_0}$ is injective.\\
(ii) Assume that $f(x)=(x-\theta_1)(x-\theta_2)(x-\theta_3)$ is
irreducible with cyclic Galois group, and with splitting field $K$ having class number $1$. Set $e:=\sqrt{D}$,
where $D$ denotes the discriminant of $f$. Let $\cal G$ denote the set
of prime factors of $eg$ in $\mathbb{Z}[t]$ that split
in
$\mathcal{O}_K[t]$. Then the specialization homomorphism $\sigma_{t_0}$ is injective, provided $t_0$ satisfies the following condition:\\
(B) For each prime factor from $\cal G$, let us choose either none
or two of its prime factors in $\mathcal{O}_K[t]$ (say
$\mathcal{P},\mathcal{Q}$). Let $h$ denote
the product of all chosen $\mathcal{P},\mathcal{Q}$. Then for each nonconstant $h$ and  each unit  $\epsilon$ of $\mathcal{O}_K$, $\epsilon h(t_0)$ is not a square from $K$.\\

\end{theorem}
\pf   Let $P\in E_g(\mathbb{Q}(t))$ be a nonzero point such that $\sigma_{t_0}(P)=O$. By the proof of Theorem \ref{main}, it is sufficient to
 prove that $P\in 2E_g(\mathbb{Q}(t))$.\\
 (i) Put $x(P)=\frac{p(t)}{q(t)^2}$ with $p,q\in \mathbb{Z}[t]$ coprime. Then
$p(p^2+apgq^2+bg^2q^4$) is a square in $\mathbb{Z}[t]$. Since $q(t_0)=0$ we see that $p(t_0)$ is a rational square.  We claim that both $p$ and
$A(t)=A=p^2+apgq^2+bg^2q^4$ are  squares in $\mathbb{Z}[t]$.
Write $p=p_0p_1^2$ where $p_0,p_1\in \mathbb{Z}[t]$ with $p_0$ square free. Then $A=p_0A_1^2$ for some $A_1\in \mathbb{Z}[t]$. Therefore $p_0$
divides $bg$ in $\mathbb{Z}[t]$. Since $p_0(t_0)$ is a non-zero rational square, we conclude, by condition (A1) of the theorem, that $p_0=1$,
 hence the claim.\\
Now we can write $A=m^2B^2$ where $m\in \mathbb{Z}$ and $B$ is a product of irreducible polynomials from $\mathbb{Z}[t]$ of positive degrees.
Also $A=(p-\theta gq^2)(p-\bar\theta gq^2)$. We claim that $p-\theta gq^2$ is a square in
$\mathcal{O}_K [t]$. By Lemma \ref{duplasto}, this and the fact that $x(P)\in \mathbb{Q}(t)^2$, imply $P\in 2E_g(\mathbb{Q}(t))$. We consider possible types of irreducible factors $H$ of $mB\in \mathbb{Z}[t]$. Let $k$ denote the multiplicity of $H$ in $B$.\\
(i1) If $H$ (constant or nonconstant) remains irreducible in $\mathcal{O}_K [t]$, then it has the same multiplicities in
$p-\theta gq^2$ and $p-\bar\theta gq^2$, say $n$. Therefore $H$ is an irreducible factor of $eg\in \mathbb{Z}[t]$. Namely, it is a factor both of
 $(\theta-\bar\theta)gq^2$ and $(\theta-\bar\theta)p$, hence it is a factor of $(\theta-\bar\theta)g$. Now we recall that
 $\theta-\bar\theta=\pm e\sqrt{d}$. We see that $2n=2k$, hence $H$ contributes in $(\theta-\bar\theta)gq^2$ with multiplicity $n=k$.  \\
(i2) If $H$ (constant or nonconstant) splits in
$\mathcal{O}_K[t]$, say $H=\pm \mathcal{P}\bar{\mathcal{P}}$, then
let $n, n'$ be the multiplicities of $\mathcal{P},
\bar{\mathcal{P}}$ in $p-\theta gq^2$. Therefore, $n+n'=2k$, hence
 $n,n'$ are even, or $H$ divides both $p-\theta gq^2$ and
$p-\bar\theta gq^2$. In the later case $H$ divides $eg$ in
$\mathbb{Z}[t]$. \\
(i3) Assume $H$ is a rational prime factor of $m$ that ramifies in
$K$, say $H=\epsilon \mathcal{P}^2$, for a prime $\mathcal{P}$ and
a unit $\epsilon$ in $\mathcal{O}_K$. Therefore the multiplicities
of $\mathcal{P}$ in $(p-\theta gq^2)\mathcal{O}_K$ and
$(p-\bar\theta gq^2)\mathcal{O}_K$ coincide (and equal to $2k$,
the multiplicity of $H$ in $m^2$).  We see that $H$
is a factor both of $(\theta-\bar\theta)gq^2$ and $(\theta-\bar\theta)p$, hence $H$ is a divisor of $eg$ in $\mathbb{Z}[t]$.\\
 By (i1), (i2) and (i3), we conclude that
 \begin{equation*}
p-\theta gq^2=\epsilon vu^2,
\end{equation*}
where $\epsilon$ is a unit in $\mathcal{O}_K$, $v$ is a factor of
$eg$ in $\mathbb{Z}[t]$ and $u\in \mathcal{O}_K[t]$. Recall that
$(p-\theta gq^2)(t_0)=p(t_0)$ is a non-zero square in
$\mathbb{Q}$, especially $\epsilon v(t_0)$ is a square in
$K$. Therefore, if $v$ is a constant polynomial, then
$\epsilon v$ is a square in $K$, hence $p-\theta gq^2$ is a
square in $\mathcal{O}_K[t]$, as we claimed. Assume now that $v$ is nonconstant.
We will show that it leads to a contradiction. Namely, in that
case, by condition (A2), and Lemma \ref{kvadrat}, we see
 $\epsilon v(t_0)$ is not a square in $K$, the contradiction. To resume, we get that $p-\theta gq^2$ is a square in
$\mathcal{O}_K [t]$.  From this, the fact that $p$ is a square
in $\mathbb{Z}[t]$ and Lemma \ref{duplasto}, we conclude that $P\in 2E_g(\mathbb{Q}(t))$, as
we claimed.\\

\bigskip

(ii) Let  $\tau$ be a generator of the Galois group of
$K/\mathbb{Q}$. We may assume that $\theta_2=\theta_1^{\tau}$,
hence $\theta_3=\theta_2^{\tau}$. Similarly as in (i) we get that
$$A=p^3+agp^2q^2+bg^2pq^4+cg^3q^6=(p-\theta_1 gq^2) (p-\theta_2 gq^2)(p-\theta_3 gq^2)$$
is a square in $\mathbb{Z}[t]$. Since $q(t_0)=0$ we see that $p(t_0)$ is a square in $K$. To prove that $\sigma_{t_0}$ is
injective it is sufficient
to prove that $p-\theta_1 gq^2$ is a square in
$\mathcal{O}_K[t]$ (by Lemma \ref{duplasto}). We can write $A=m^2B^2$ as in (i). Let $H$ be an irreducible factor of $mB\in \mathbb{Z}[t]$ of multiplicity $k$.\\
(ii1) If $H$  (constant or nonconstant)  remains irreducible in
$\mathcal{O}_K [t]$, then it has the same multiplicity, say $n$ in
each $p-\theta_i gq^2,\ i=1,2,3$. Therefore the multiplicity is
even. On the other side $H^n$ divides $g$, hence, if $H$ is nonconstant, then $n=0$.\\
(ii2) If $H$  (constant or nonconstant) splits in
$\mathcal{O}_K[t]$, say
$H=\pm\mathcal{P}\mathcal{P}^{\tau}\mathcal{P}^{\tau^2}$, then let
$n,n',n''$ be the corresponding multiplicities of $\mathcal{P},
\mathcal{P}^{\tau},\mathcal{P}^{\tau^2}$ in $p-\theta_1 gq^2$.
Therefore $n+n'+n''=2k$, hence the multiplicities are even, or two
of them are odd. In the later case we get that $H$ divides $eg$ in $\mathbb{Z}[t]$. \\
(ii3) If $H$ is is constant and ramifies in $K$, then
$H=\epsilon\mathcal{P}^3$. Let $n$ be the multiplicity of
$\mathcal{P}$ in $p-\theta_1 gq^2$. We get $n=2k$.\\
 By (ii1), (ii2) and (ii3) we see that there exist unit $\epsilon\in
 \mathcal{O}_K$, $h\in \mathcal{O}_K[t]$ as in (B), and
and $w\in \mathcal{O}_K[t]$ such that
$$p-\theta_1 gq^2=\epsilon\cdot hw^2.$$
The assumption that $h$ is a nonconstant polynomial leads to a contradiction (recall that $p(t_0)$ is a square in $K$, $q(t_0)=0$
and $\epsilon h(t_0)$ is not a square in $K$). Therefore, $h$ is constant, which implies that
 $p-\theta_1 gq^2$ is a square in $\mathcal{O}_K[t]$, as we claimed.
\qed

\section{An example by Mestre}
\label{section5}

In the former sections we presented the main results of this paper. In the next section we will show some ways of applying the results to concrete examples.  We will use the results for calculating the rank and proving that a set of points are free generators of an elliptic curve over the field of rational functions in one variable, by choosing a particular specialization which we will know is injective by Theorem \ref{main} or Theorem \ref{tmtwist}. Before that, in this section we need to mention a few things.

In \cite{Me} and \cite[Theorem 3.7]{R-S1} (see also \cite[Theorem 3]{Ste-T}), the following has been shown for a family of twists of the elliptic curve over $\mathbb Q(u)$, which we will observe in the next section.

\begin{example}\label{exam}
Let $a,b\in\mathbb Q$ such that $ab\ne 0$, let
$$g(u)=g^{a,b}(u)=-ab\cdot (u^2+1)\cdot (b^2(u^4+u^2+1)^3+a^3u^4(u^2+1)^2)$$
 and let $E^{a,b}$ be the elliptic curve over $\mathbb Q$ given by the equation
 $$E^{a,b}:y^2=x^3+ax+b.$$
 Then $E^{a,b}_g:y^2=x^3+ag(u)^2x+g(u)^3$  has rank at least $2$ over $\mathbb Q(u)$, with two independent points $P_g^{a,b}$ and $Q_g^{a,b}$
 with  coordinates
 \begin{equation}\label{P}
 P_g^{a,b}=\left(-\frac ba\frac{(u^2+u+1)(u^2-u+1)}{(u^2+1)}\cdot g(u),\frac 1{a^2(u^2+1)^2}\cdot g(u)^2\right)
 \end{equation}
 \begin{equation}\label{Q}
 Q_g^{a,b}=\left(-\frac ba\frac{(u^2+u+1)(u^2-u+1)}{u^2(u^2+1)}\cdot g(u),\frac 1{a^2u^3(u^2+1)^2}\cdot g(u)^2\right).
 \end{equation}
\end{example}
\medskip

It is easy to see that
{\small{
$$g^{a,b}(u)=-ab(u^2+1)(b^2u^{12}+3b^2u^{10}+(6b^2+a^3)u^8+(7b^2+2a^3)u^6+(6b^2+a^3)u^4+3b^2u^2+b^2).$$}}
Put,
$$f(x)=x^3+ax+b=(x-\theta_1)(x-\theta_2)(x-\theta_3).$$
 Let $K$ denote the splitting field of $f$. By using the point \eqref{P} we conclude that $g=g^{a,b}$ factors over $K$ into
 {\small{
$$-a(u^2+1)( bu^4+(b+a\theta_1)u^2+b+a\theta_1)(bu^4+(b+a\theta_2)u^2+b+a\theta_2)(bu^4+(b+a\theta_3)u^2+b+a\theta_3).$$}}
Namely,
{\scriptsize{
$$g(u)=a^4(u^2+1)^4\left(-\frac ba\frac {u^4+u^2+1}{u^2+1}- \theta_1\right)\left(-\frac ba\frac {u^4+u^2+1}{u^2+1}-\theta_2\right)\left(-\frac ba\frac {u^4+u^2+1}{u^2+1}-\theta_3\right).$$
}}
 Therefore, if $K$ is a cubic field (with cyclic Galois group), then  the set $\cal G$ from Theorem \ref{tmtwist} (ii)  is nonempty. Moreover, it contains non-constant elements, especially the implementation of the criterion is nontrivial. Note also that, in this case,  $g^{a,b}$ factors over $\mathbb{Q}$ as
$$g^{a,b}(u)=-ab(u^2+1)(bu^6+\frac{3b-e}{2}u^4+\frac{3b+e}{2}u^2-b)(bu^6-\frac{3b+e}{2}u^4-\frac{3b-e}{2}u^2-b),$$
where $e=\sqrt{D}$ as in Theorem \ref{tmtwist}.

Further we will use the following lemmas, which are easy to prove.

\begin{lemma}\label{submain}
Let $\sigma:G\rightarrow H$ be a homomorphism of finitely generated abelian groups.
If  $\sigma$ is an injection and  the rank of $G$ is greater or equal to the rank of $H$, then  $G$ and $H$ have the same rank.
\end{lemma}

\begin{lemma}\label{inj}
Let $\sigma:G\rightarrow H$ be a homomorphism of finitely generated abelian groups which is injective on the torsion subgroup. Let $P_1,\ldots ,P_r$ be the free generators of $G$. Then $\sigma(P_1),\ldots ,\sigma(P_r)$ are independent over $\mathbb Z$ if and only if $\sigma$ is injective.
\end{lemma}

For a point $T$ on $E^{a,b}_g$, we will denote by $T'$ the corresponding point on the model $g(u)y^2=f(x)$ of $E^{a,b}_g$. It can be checked that\\

$(P_g^{a,b})'=(-\frac{b(u^2+u+1)(u^2-u+1)}{a(u^2+1)}, \frac{1}{a^2(u^2+1)^2})$,\\

$(Q_g^{a,b})'=(-\frac{b(u^2+u+1)(u^2-u+1)}{au^2(u^2+1)},
\frac{1}{a^2u^3(u^2+1)^2})$,\\

$(P_g^{a,b}+Q_g^{a,b})'=(-\frac{b^2u^8-2b^2u^7+(4b^2+a^3)u^6-4b^2u^5+(5b^2+2a^3)u^4-4b^2u^3+(4b^2+a^3)u^2-2b^2u+b^2}
{ab(u^2+1)(u+1)^2(u^2-u+1)^2},$\\

$\frac{b^2u^6-(3b^2+a^3)u^5+6b^2u^4-(7b^2+2a^3)u^3+6b^2u^2-(3b^2+a^3)u+b^2}{a^2b^2(u^2+1)^2(u+1)^3(u^2-u+1)^3})$,\\

$(P_g^{a,b}-Q_g^{a,b})'=(-\frac{b^2u^8+2b^2u^7+(4b^2+a^3)u^6+4b^2u^5+(5b^2+2a^3)u^4+4b^2u^3+(4b^2+a^3)u^2+2b^2u+b^2}
{ab(u^2+1)(u-1)^2(u^2+u+1)^2},$\\

$-\frac{b^2u^6+(3b^2+a^3)u^5+6b^2u^4+(7b^2+2a^3)u^3+6b^2u^2+(3b^2+a^3)u+b^2}{a^2b^2(u^2+1)^2(u-1)^3(u^2+u+1)^3}).$\\

\begin{lemma}
 Let $E:y^2=x^3+ax^2+bx+c$ be an elliptic curve over a number field
$k$, and let $g$ be a nonconstant square free polynomial from $k[u]$.
Let $E_g:y^2=x^3+ag(u)x^2+bg(u)^2x+cg(u)^3$ be the quadratic twist of $E$ by $g$.
For $T\in E_g(k(u))$ let us define $h_0(u):=\deg \frac{x(u)}{g(u)}$, where $x(u)$ denotes the $x$-coordinate of $T$ as a rational function in $u$.\\
Then $\frac{1}{2}h_0$ is the canonical height on $E_g$.
\end{lemma}
\pf \cite[Theorem 1.]{G-L} \qed  \medskip

Let $\langle\ ,\ \rangle$ denote the canonical bilinear form on $E_g(k(u))\times E_g(k(u))$. Then
$$\langle T,S\rangle=\frac{1}{2}(h_0(T+S)-h_0(T)-h_0(S)),\ {\rm for\ each}\ T,S\in E_g(k(u)),$$
especially
$$h_0(T)=\langle T,T\rangle,\ {\rm for\ each}\ T.$$
By \cite{Ste-T}, Section 4, Corollary 1, and
\cite{R-S1}, Remark 2.12, we see that the rank of $E^{a,b}_g$ over $\mathbb{Q}(u)$ is $\leq 6$.

In the following lemma we prove that $P^{a,b}_g,Q^{a,b}_g$ are free generators of $E^{a,b}_g(\mathbb{Q}(u))$ provided the rank is two.

\begin{lemma}\label{GL}
Let $P^{a,b}_g,Q^{a,b}_g$ be points from Example \ref{exam}. Assume that the rank of $E^{a,b}_g(\mathbb{Q}(u))$ is equal to $2$. Then $P^{a,b}_g,Q^{a,b}_g$
are free generators of $E^{a,b}_g(\mathbb{Q}(u))$.
\end{lemma}
\pf We will repeat the argument from \cite[Section 3]{G-L}.
Let $T\in E_g(\mathbb{Q}(u))$ be an arbitrary nonzero point. We will show that $T$ is a $\mathbb{Z}$-linear combination of $P^{a,b}_g,Q^{a,b}_g$ and torsion points. Since
the rank is supposed to be two, the points $P^{a,b}_g,Q^{a,b}_g,T$ are dependent, so
 there is a nontrivial relation
\begin{equation*}
kT=mP^{a,b}_g+nQ^{a,b}_g,\ m,n,k\in \mathbb{Z}.
\end{equation*}
We may assume that $k\geq 0$, hence $k>0$ (contrary, the relation
is trivial). By consecutive adding or subtracting $kP^{a,b}_g$ or $kQ^{a,b}_g$, it
leads to
\begin{equation*}
kT_1=m'P^{a,b}_g+n'Q^{a,b}_g,\ {\rm with}\ -\frac{k}{2}\leq m',n'\leq \frac{k}{2}
\end{equation*}
(for example, if we start with $5T=13P^{a,b}_g-4Q^{a,b}_g$ we get $5T_1=-2P^{a,b}_g+Q^{a,b}_g$
where $T_1=T-3P^{a,b}_g+Q^{a,b}_g$). We get
\begin{equation*}
k^2h_0(T_1)=m'^2h_0(P^{a,b}_g)+n'^2h_0(Q^{a,b}_g)+2m'n'\langle P^{a,b}_g,Q^{a,b}_g\rangle,
\end{equation*}
which provides an upper bound for $h_0(T_1)$:
\begin{equation}\label{h0T1}
h_0(T_1)\leq \frac{h_0(P^{a,b}_g)+h_0(Q^{a,b}_g)+2|\langle P^{a,b}_g,Q^{a,b}_g\rangle|}{4}.
\end{equation}
In our case $h_0(P^{a,b}_g)=h_0(Q^{a,b}_g)=4$ (since the expressions for $x(P^{a,b}_g)$ and $x(Q^{a,b}_g)$ are reduced to lowest terms for each $a,b$). Also $h_0(P^{a,b}_g+Q^{a,b}_g)\leq 8$ and $h_0(P^{a,b}_g-Q^{a,b}_g)\leq 8$, for each $a,b$. By the parallelogram law (recall that $\frac{1}{2}h_0$ is the canonical height), we conclude that $h_0(P^{a,b}_g+Q^{a,b}_g)=h_0(P^{a,b}_g-Q^{a,b}_g)=8$,
hence $ \langle P^{a,b}_g,Q^{a,b}_g\rangle=0$. By \eqref{h0T1} we get $h_0(T_1)\leq 2$. We claim that $h_0(T_1)=1$ or $h_0(T_1)=2$ is impossible.  Contrary,
 $x=x(T_1)/g(u)=\frac{\alpha(u)}{\beta(u)}$, where $\alpha,\beta$ are nonzero polynomials over $\mathbb{Q}$ of degree at most $2$ and at least one of them is non-constant. Plugging in $g(u)y^2=x^3+ax+b$, we get that there is a nonzero polynomial $w$ over $\mathbb{Q}$ of degree at most $6$, such that $w(u)\beta(u)g(u)$ is a square in $\mathbb{Q}[u]$. It is impossible since $g(u)$ is squarefree with degree $14$. Therefore $h_0(T_1)=0$, hence $T_1$ is torsion, i.e. $T$ is a $\mathbb{Z}$-linear combination of $P^{a,b}_g,Q^{a,b}_g$ and torsion points. The lemma is proved. \qed
 \medskip

\section{Application of our method to the Mestre example}
\label{section6}

In this last section we will use the results obtained in the former
sections to  concrete examples, by using Theorem \ref{main} for $K=\mathbb Q$ and Theorem \ref{mainh} for $K$ an algebraic number field of class number one and two, and Theorem \ref{tmtwist} for $K$ an algebraic number field of class number one.

We will calculate the rank and prove
that a given set of points are free generators of some elliptic
curves over the field of rational functions in one variable over
$\mathbb Q$, i.e $\mathbb Q(u)$. For this we will use a concrete
family (of twists) of elliptic curves over $\mathbb Q(u)$
mentioned in Example \ref{exam}.

Here we will show that for certain values of $a,b\in\mathbb Z$,
the elliptic curves  $E^{a,b}_g$ over $\mathbb Q(u)$, coming from
Example \ref{exam}, have rank two and free generators the two
points mentioned, by using Theorem \ref{main}, Theorem \ref{mainh} or Theorem
\ref{tmtwist} to pick a particular specialization which we will
know is injective. We  observe integer values $a,b$,  depending
on the Galois group of the polynomial $f(x)=x^3+ax+b$. We will
observe only $f$ having  splitting fields  $K$ with  class number
one and two. The case $K=\mathbb Q$  (i.e. in cases when the right side of
the equation of the curve $E^{a,b}_g$ factors over $\mathbb Q(u)$
into linear factors), we will call  {\it the rational case}.  The
case when $K$ is  a quadratic field (i.e. when the right side of the equation
of the curve $E^{a,b}_g$ over $\mathbb Q(u)$ factors into exactly
two factors), we will call  {\it the quadratic case}.  For $f$
with  Galois group the cyclic  group of order three  or the
symmetric group ${\bf S_3}$ (in cases when the right side of the equation of the
curve $E^{a,b}_g$ is irreducible over $\mathbb Q(u)$), we will
talk about {\it the  cubic cyclic case}  or {\it the symmetric
case}.

We use Theorem \ref{main} in the rational and Theorem \ref{mainh} in the symmetric case of class number one. We
use Theorem \ref{tmtwist} (i) in the quadratic case, for
$K=\mathbb Q(\sqrt{d})$ of class number one such as
$$d=-163,-67,-43,-19,-11,-7,-3,-2,-1,3,6,7,11,13,14,17,19,$$
(among them are all imaginary quadratic fields of class number
one), and Theorem \ref{tmtwist} (ii) in some cyclic cases of class number one.
Finally at the end is presented the usage of Theorem \ref{mainh} for a few  quadratic case when the splitting field $K=\mathbb Q(\sqrt {-5})$ is of class number two.

Calculations in this section were performed using a variety of packages: GP/Pari \cite{Pari}, MAGMA
\cite{MAGMA}, mwrank \cite{mwrank}.

So the following are the cases when the splitting field of $E^{a,b}$ is of class number one.

\begin{theorem} \label{corall}
Let $a,b$ be given in one of the following ways:
\begin{itemize}
\item[a)]
$$a=r_1r_2-(r_1+r_2)^2,\ b=r_1r_2(r_1+r_2)$$
where $1\leq r_1< r_2\leq 15.$
\item[b)]
$$a=-3 r_1^2-d\cdot r_2^2,\ b=2 r_1\cdot (r_1^2-d\cdot r_2^2),$$
where $1\leq r_1,r_2\leq 3,$
{\small{
$$d=-163,-67,-43,-19,-11,-7,-3,-2,-1,3,6,7,11,13,14,17,19,$$
}}
and for $d=-3$ we take $r_1\ne r_2$ and $(r_1,r_2)\ne (1,3)$.
\item[c)]
 $a,b$ are the ones in the third  main column of the Table.
\item[d)] $a,b$ are the ones in the last  main column of the
Table.
\end{itemize}
Let $E^{a,b}$ be the elliptic curve given with the equation over
$\mathbb Q$
$$E^{a,b}:y^2=x^3+ax+b.$$
Then the elliptic curve $E^{a,b}_g$   has rank
two over $\mathbb Q(u)$, with free generators the two points $P^{a,b}_g$ and
$Q^{a,b}_g$.

{\rm{
{\tiny{
\begin{tabular}{|||c|c|c|||r|c|c|c|||c|c|c|||c|c|c|||}
\multicolumn{3}{c}{Rational}&\multicolumn{4}{c}{Quadratic}&\multicolumn{3}{c}{Cyclic cubic}&\multicolumn{3}{c}{Symmetric}\\
\hline \hline \hline
$a$&$b$&$u_0$&$d$&$a$&$b$&$u_0$&$a$&$b$&$u_0$&$a$&$b$&$u_0$\\
\hline
-7&6&14&-163&160&328&8&-189&189&5&1&1&3\\
-13&12&20&-163&649&1306&32&-837&3348&4&4&4&2\\
-21&20&64&-67&64&136&5&-5697&74061&3&8&8&2\\
-31&30&14&-67&265&538&5&-8289&132624&3&9&9&2\\
-43&42&10&-43&40&88&12&-11367&215973&2&22&44&2\\
-57&56&18&-43&169&346&9&-46521&1860840&3&36&144&3\\
-73&72&15&-19&16&40&9&-1809&9045&3&63&441&4\\
-91&90&34&-19&73&154&6&-5373&59103&6&&&\\
-111&110&9&-11&8&24&15&-39069&1367415&4&&&\\
-133&132&17&-11&41&90&6&-3&-1&3&&&\\
-157&156&11&-7&4&16&2&-333&999&4&&&\\
-183&182&10&-7&25&58&8&-1647&1647&10&&&\\
-211&210&6&-3&9&26&6&-4887&34209&5&&&\\
-241&240&19&-3&-9&28&6&-10071&130923&4&&&\\
-19&30&9&-2&-1&6&4&-37287&1155897&3&&&\\
-28&48&14&-2&5&18&19&-50247&1859139&4&&&\\
-39&70&13&-1&-2&4&9&-513&-3591&2&&&\\
-52&96&20&-1&1&10&21&-2943&5886&3&&&\\
-67&126&10&3&-6&-4&3&-9747&107217&2&&&\\
-84&160&64&3&-15&-22&5&-16713&284121&5&&&\\
-103&198&12&6&-9&-10&10&-42633&1364256&3&&&\\
-124&240&14&6&-27&-46&5&-9&-9&3&&&\\
-147&286&22&7&-10&-12&8&-441&1323&4&&&\\
-172&336&10&7&-31&-54&5&-999&-10989&3&&&\\
-199&390&29&11&-14&-20&10&-12339&123390&3&&&\\
-228&448&18&11&-47&-86&5&-29511&649242&3&&&\\
-259&510&39&13&-16&-24&2&-35019&875475&4&&&\\
-37&84&5&13&-55&-102&12&-47493&1472283&4&&&\\
-49&120&28&14&-17&-26&6&-1323&-17199&4&&&\\
-63&162&14&14&-59&-110&5&-9099&45495&7&&&\\
-79&210&9&17&-20&-32&5&-15579&171369&4&&&\\
-97&264&17&17&-71&-134&2&-2133&-36261&3&&&\\
-117&324&20&19&-22&-36&6&-8937&8937&2&&&\\
-139&390&16&19&-79&-150&9&-23193&301509&2&&&\\
\hline \hline \hline
\end{tabular}
}}
}}
\begin{center}
{\bf Table 4.2} List of  values $a,b$ with corresponding $u_0$ for which the specialization $\sigma_{u_0}$ is injective and  $E_{g}^{a,b}(u_0)(\mathbb Q)$ has rank two. This presents just a portion of values $a,b$ treated in the theorem.
\end{center}
\end{theorem}

\pf

Let $K$ be the splitting field of $f(x)=x^3+ax+b$, in all cases from the claim of the theorem $K$ is of class number one.
We thus observe four cases:
\begin{itemize}
\item the {\bf rational case} is presented in a):\\
for integers
$$a=r_1r_2-(r_1+r_2)^2,\ b=r_1r_2(r_1+r_2)$$
 as in the statement of the theorem  we have
$$E^{a,b}:y^2=(x-r_1)(x-r_2)(x+r_1+r_2)=x^3+ax+b,$$
and  the corresponding splitting field is $K=\mathbb Q$.
\item the {\bf quadratic case} is presented in b):\\
for integers
$$a=-3 r_1^2-d\cdot r_2^2,\ b=2 r_1\cdot (r_1^2-d\cdot r_2^2),$$
as in the statement of the theorem we have
$$E^{a,b}:y^2=(x+2r_1)(x^2-2r_1x+r_1^2-d\cdot r_2^2)=x^3+ax+b$$
and  the corresponding splitting field is $K=\mathbb Q(\sqrt d)$.
\item the {\bf cyclic cubic case} is presented c), \item the {\bf
symmetric case} is presented d).
\end{itemize}

We proved that for the values $a,b$ in the statement of the
Theorem there exists an integer $u_0$  for which the conditions
of the Lemma \ref{submain} are satisfied for the specialization
homomorphism $\sigma_{u_0}$. In other
words:
\begin{itemize}
\item we choose a specialization $\sigma_{u_0}$ from
specializations which we know are  injections by using  Theorem
\ref{main}, Theorem
\ref{mainh} or Theorem \ref{tmtwist}, in the following way. We take the splitting field $K$ of $f(x)=x^3+ax+b$, in all cases treated it is of class number one so we can apply Theorem
\ref{main}, Theorem
\ref{mainh} or Theorem \ref{tmtwist}. In
the rational case we apply Theorem \ref{main} for $K=\mathbb Q$,
in the quadratic   case we apply Theorem \ref{tmtwist} (i), in the cubic cyclic Theorem \ref{tmtwist} (ii),
and in the symmetric case we apply Theorem \ref{mainh} for $K$. Note that in the quadratic and cyclic cubic
case we also could have applied Theorem \ref{main} with $K$ the
corresponding splitting field.
\item  we check if the  root number is one and then if the rank of the specialized curve $E^{a,b}_g(u_0)$ is 2 over $\mathbb Q$
(which is calculated with {\tt mwrank} for the rational and
quadratic case, with  Magma's command \\ {\tt MordellWeilShaInformation} for the cubic cyclic and the symmetric
case)
\end{itemize}
Thus for these values of $a,b$ and the chosen specialization
homomorphism $\sigma_{u_0}:E^{a,b}_g(\mathbb Q(u))\rightarrow E^{a,b}_g(u_0)(\mathbb Q)$ given by $u\mapsto u_0$, Lemma \ref{submain} is
satisfied, since we know that the rank of $E^{a,b}_g$ over $\mathbb
Q(u)$ is at least two (since we have two independent points). So we conclude that it is actually exactly
two. Now  Lemma \ref{GL}  is applied to conclude that $E^{a,b}_g$
over $\mathbb Q(u)$ has actually free generators the two points
$P^{a,b}_g,$ $Q^{a,b}_g$. In the table above are listed some chosen
values $u_0$, depending on $a,b$.

 For the symmetric case we looked at $\sigma_{u_0}=\sigma_{u_0}^K\vert_{E^{a,b}_g(\mathbb
Q(u))}:E^{a,b}_g(\mathbb Q(u))\rightarrow E^{a,b}_g(u_0)(\mathbb Q)$, where $\sigma_{u_0}^K:E^{a,b}_g(K(u))\rightarrow E^{a,b}_g(u_0)(K)$ is proved to be injective by Theorem \ref{mainh}.

\medskip

a) {\bf{ We present the proof for one rational case}}:  $a=-7,\ b=6.$\\
This corresponds to  $r_1=1,\ r_2=2$. Thus the corresponding
elliptic curve $E^{-7,6}$ over $\mathbb Q$ is
$$E^{-7,6}:y^2=(x-1)(x-2)(x+3)=x^3-7x+6.$$ Then
$$g^{-7,6}(u)=1512u^{14}+6048u^{12}-798u^{10}-23562u^8-23562u^6-798u^4+6048u^2+1512=$$
$$=2\cdot 3\cdot 7\cdot(u^2 - 2)(u^2 + 1)(u^2 + 3)(2u^2 - 1)(2u^2 + 3)(3u^2 + 1)(3u^2 + 2).$$
Thus $$e_1(u)=1\cdot g^{-7,6}(u),\ \ \ e_2(u)=2\cdot g^{-7,6}(u),\
\ \ e_3(u)=-3\cdot g^{-7,6}(u),$$ and so

\noindent
 {\small{
$\mbox{rad}_{\mathbb Z[u]}(e_1(u)-e_2(u))\cdot(e_1(u)-e_3(u))=$
$$2\cdot 3\cdot 7\cdot (u^2 + 1)(u^2 - 2)(u^2 + 3)(2u^2 - 1) (2u^2 + 3)(3u^2 + 1)(3u^2 + 2),$$
$\mbox{rad}_{\mathbb Z[u]}(e_2(u)-e_1(u))\cdot(e_2(u)-e_3(u))=\mbox{rad}_{\mathbb Z[u]}(e_3(u)-e_1(u))\cdot(e_3(u)-e_2(u))=$
$$=2\cdot 3 \cdot 5\cdot 7\cdot (u^2 + 1)(u^2 - 2)(u^2 + 3)(2u^2 - 1) (2u^2 + 3)(3u^2 + 1)(3u^2 + 2).$$

For example, to obtain all nonconstant square-free divisors of $e_1(u)-e_2(u))\cdot(e_1(u)-e_3(u)$ in $\mathbb Z[u]$ we look at all nonconstant
{\small{
$$(-1)^{i_1}\cdot 2^{i_2}\cdot 3^{i_3}\cdot 7^{i_4}\cdot (u^2 + 1)^{i_5}(u^2 - 2)^{i_6}(u^2 + 3)^{i_7}(2u^2 - 1)^{i_8} (2u^2 + 3)^{i_9}(3u^2 + 1)^{i_{10}}(3u^2 + 2)^{i_{11}},$$}}where $i_k\in\{0,1\}$, for $k=1,\ldots 11.$

If we choose $u_0=14$ we have
\begin{itemize}
\item  it is easy to see  that if $h(u)$ is a nonconstant divisor of rad$_{\mathbb Z[u]}((e_1(u)-e_2(u))\cdot(e_1(u)-e_3(u)))$ or rad$_{\mathbb Z[u]}((e_2(u)-e_1(u))\cdot(e_2(u)-e_3(u)))$ or rad$_{\mathbb Z[u]}((e_3(u)-e_1(u))\cdot(e_3(u)-e_2(u)))$ in $\mathbb Z[u]$, then $h(14)$ is not a square in $\mathbb Q.$ The divisors of $g$ are checked only once. \\
Thus Theorem \ref{main} is satisfied, so we conclude that the
specialization homomorphism $\sigma_{14}:E^{-7,6}_g(\mathbb
Q(u))\rightarrow E^{-7,6}_g(14)(\mathbb Q)$ is an injection.
 \item
the elliptic curve $E^{-7,6}_g(14)$ over $\mathbb Q$ is given by the equation {\small{
$$y^2=x^3-2057410750080462983177474912957475480000x+$$
$$+30233310019074218054503104857297537715567578648874672000000,$$}}
\item {\tt {mwrank}} \cite{mwrank} showed that  $E^{-7,6}_g(14)$ has rank $2$ over $\mathbb Q$.
\end{itemize}
Thus for   $(a,b)=(-7,6)$ and  the chosen specialization homomorphism $\sigma_{14}$, Lemma \ref{submain} is satisfied, so we conclude
that the rank of $E^{-7,6}_g$ over $\mathbb Q(u)$ is two. Now by
applying Lemma \ref{GL}  we conclude that $E^{-7,6}_g$ over
$\mathbb Q(u)$ has rank two and free generators
$P^{-7,6}_g,Q^{-7,6}_g$.

\medskip

b) {\bf {We present the proof for one quadratic case:}} $a=1, \ b=10.$\\
It corresponds to  $r_1=1,\
r_2=2,\ d=-1$. Thus
$$E^{1,10}:y^2=(x+2)(x^2-2x+5)=x^3+x+10.$$
We thus deal with the splitting field the quadratic number field
$K=\mathbb Q(\sqrt{-1})$ with class number one. We take $u_0=21.$
\begin{itemize}
\item we easily check that   each non-constant divisor $h(u)$ of
rad$((9c^2-qd^2)g(u))$ or rad$(2dg(u))$ in $\mathbb Z[u]$, the
value $d'\cdot h(21)$  is not a square in $\mathbb Q$, where $d'=1$ or $d'$ each square free divisor of $4$ respectively. Which is the
condition (A1) and (A2) from Theorem \ref{tmtwist} (i). Thus the
specialization $\sigma_{21}$ given by $u_0=21$ is a monomorphism.
\item the specialized elliptic curve $E^{1,10}_g(21)$ over $\mathbb Q$ is given by
the equation

$$E^{1,10}_g:y^2=x^3+10715763113679070635989488456722194251878400x$$
{\small{
$$-350779864964306170166930397804810220833346107858355456012779520000.$$
}}
\item   {\tt{mwrank}} \cite{mwrank} showed the specialized elliptic
curve $E^{1,10}_g(21)$ over $\mathbb Q$ has rank 2 .
\end{itemize}

Thus for these values of $a,b$ and the chosen specialization
homomorphism given by $u_0=21$ Lemma \ref{submain} is satisfied, so
we conclude that the rank of $E^{1,10}_g$ over $\mathbb Q(u)$ is
two. Now applying Lemma \ref{GL}  we conclude that $E^{1,10}_g$
over $\mathbb Q(u)$ has rank two and free generators
$P^{1,10}_g,Q^{1,10}_g$.

\medskip

c) {\bf{ We present the proof for one cubic cyclic case}}: {\small {for $a=-1647,\ b=1647.$}}\\
We have
$$E^{-1647,1647}:y^2=x^3-1647x+1647,$$
 $e=3^7\cdot 61$ (see Theorem \ref{tmtwist} (ii)). So we observed the cyclic cubic field
$K=\mathbb{Q}(q)$ with class number one where $q$ is a root of the polynomial
$x^3-1647x+1647$, which has two fundamental units (which were found using the command {\tt {bnfinit}} in Pari)
$$\epsilon_1=\frac 1{81} q^2 - \frac {13}{27}q + \frac 49,$$
$$\epsilon_2=\frac 1{81}q^2+ \frac{14}{27}q - \frac 59.$$

The corresponding function  is
{\small{
$$g^{-1647,1647}(u)=7358247586881u^{14} + 29432990347524u^{12} - 12052809547311078u^{10}$$
$$ - 36261444108149568u^8 - 36261444108149568u^6 - 12052809547311078u^4$$
$$ + 29432990347524u^2 + 7358247586881.$$
}}
There are five irreducible factors of $e\cdot g^{-1647,1647}$ in
$\mathbb Z[u]$, these are
$$3, 61, u^2 + 1, u^6 - 39u^4 - 42u^2 - 1, u^6 + 42u^4 + 39u^2 - 1.$$
We check which are their prime factors in $\mathcal O_K[u]$ using Pari:
{\tiny{
\begin{itemize}
\item $3:\ \ \ \ \ \  -\frac 1{243}q^2 - \frac{14}{81}q +
\frac{14}{27}, \frac 1{243}q^2 + \frac{14}{81}q + \frac{13}{27},
\frac 2{243}q^2 +\frac{ 28}{81}q - \frac 1{27}$ \item
$61:\ \ \ \ \ \ \frac 1{81}q^2 -\frac 4{27}q - \frac{122}9$ \item
$u^2+1:\ \ \ \  u^2+1$
\item $u^6 - 39u^4 - 42u^2 - 1:\ \
u^2 -\frac 1{81}q^2 -\frac  {14}{27}q + \frac 59, u^2 -\frac
1{81}q^2 + \frac{13}{27}q + \frac 59, u^2 + \frac 2{81}q^2 + \frac
1{27}q - \frac {361}9$
\item $u^6 + 42u^4 + 39u^2 - 1:\ \ \
u^2 -\frac 2{81}q^2 - \frac 1{27}q + \frac {370}9,u^2 +\frac
1{81}q^2 - \frac {13}{27}q + \frac 49, u^2+ \frac 1{81}q^2 +
\frac{14}{27}q + \frac 49$.
\end{itemize}
}}

So to check that Theorem \ref{tmtwist} (ii) is satisfied, we have
to show that
$$(-1)^{i_0}\cdot \epsilon_1^{i_1}\cdot \epsilon_2^{i_2}\cdot\prod'_{\cal P,\cal Q}\cal P\cdot \cal Q, $$
when specialized to $u_0$ is not a square in $K$, for $i_0,i_1,i_2=0,1$ and the product $\prod'_{\cal P,\cal Q}\cal P\cdot \cal Q, $ is taken as in Theorem \ref{tmtwist} (ii).

For example, for $u_0=0$ we conclude that Theorem \ref{tmtwist}
(ii) isn't satisfied. Since if we take for $h(u)$ no prime factors
from 3, 61, $u^2+1$ and $u^6 - 39u^4 - 42u^2 - 1$, and we take
from $u^6 + 42u^4 + 39u^2 - 1$ the two prime factors
$$u^2 -\frac 2{81}q^2 - \frac 1{27}q + \frac{370}9,u^2 +\frac 1{81}q^2 - \frac{13}{27}q + \frac 49,$$
 further if we take for the unit $\epsilon$ the unit part $$-\epsilon_1\cdot\epsilon_2=-\left(\frac 1{81}\cdot q^2 - \frac {13}{27}q + \frac 49\right)\left(\frac 1{81}q^2+ \frac{14}{27}q - \frac 59\right)=1-q.$$ So $\epsilon h(u)$ is equal to the invertible $1-q$ multiplied by the chosen primes, then
 $$\epsilon h(u)=(1-q)\left(u^2 -\frac 2{81}q^2 - \frac 1{27}q + \frac{370}9\right)\left(u^2 +\frac 1{81}q^2 - \frac{13}{27}q + \frac 49\right).$$
 Thus $\epsilon h$ specialized to $u_0=0$ we get a square in $K$, more precisely
 $$\epsilon h(0)=\left(-\frac 1{81}q^2 + \frac{13}{27}q -\frac 49\right)^2.$$

We list what other candidates  $\sigma_{u_0}$ give:

\begin{itemize}
\item $u_0=0,1,2,3,4,5$ doesn't satisfy Theorem \ref{tmtwist}
(ii), \item  $u_0=6,7$ satisfies Theorem \ref{tmtwist} (ii), but
doesn't have root number one, \item $u_0=8$ does satisfy Theorem
\ref{tmtwist} (ii) and has root number one, but Magma  returns that the rank is  $\geq 3$,
\item$u_0=9$ satisfies the Theorem, but doesn't have root number
one.
\end{itemize}
While $u_0=10$ satisfies Theorem \ref{tmtwist} (ii),   and
$E^{-1647,1647}_g(10)$ over $\mathbb Q$ has rank two, we conclude
by Lemma \ref{submain} that  $E^{-1647,1647}_g$ over $\mathbb
Q(u)$ has rank two. By Lemma \ref{GL} we conclude that
$E^{-1647,1647}_g$ has free generators
$P^{-1647,1647}_g,Q^{-1647,1647}_g$.

d) {\bf{ We present the proof for one symmetric case}}: {\small {for $a=1,\ b=1.$}}\\

Thus in the notation from Example \ref{exam}
$$a=1,\ \ \ b=1,$$
and the corresponding elliptic curve $E^{1,1}$ over $\mathbb Q$ is
$$E^{1,1}:y^2=x^3+x+1.$$ Then the splitting field $K$  of $x^3+x+1$ is of class number one, generated by  the algebraic number $q$ defined as a root of  $x^6 + 78x^4 + 324x^3 + 1521x^2 + 12636x + 64219$. Thus $K=\mathbb Q(q)$. The two  fundamental units of $K$ are
$$\frac 4{245805}q^5 - \frac{169}{737415}q^4 + \frac{52}{49161}q^3 - \frac{7097}{737415}q^2 - \frac{8728}{49161}q - \frac{13156}{737415}, $$
$$\frac 2{49161}q^5 - \frac{169}{294966}q^4 + \frac{130}{49161}q^3 - \frac{7097}{294966}q^2 + \frac{5521}{98322}q - \frac{6578}{147483}.$$
and further,
$$g^{1,1}(u)=-u^{14} - 4u^{12} - 10u^{10} - 16u^8 - 16u^6 - 10u^4 - 4u^2 - 1.$$

We have
 $x^3+x+1=(x-e_1(u))(x-e_2(u))(x-e_3(u)),$
where
{\scriptsize{
$$e_1(u)=-\left(-\frac 2{35115}q^5 + \frac{169}{210690}q^4 - \frac{26}{7023}q^3 + \frac{7097}{210690}q^2 + \frac{1705}{14046}q + \frac{6578}{105345}\right)\cdot g^{1,1}(u),$$
$$e_2(u)=-\left(\frac 4{245805}q^5 - \frac{169}{737415}q^4 + \frac{52}{49161}q^3 - \frac{7097}{737415}q^2 - \frac{8728}{49161}q - \frac{13156}{737415}\right)\cdot g^{1,1}(u),$$
$$e_3(u)=-\left(\frac 2{49161}q^5 - \frac{169}{294966}q^4 + \frac{130}{49161}q^3 - \frac{7097}{294966}q^2 + \frac{5521}{98322}q - \frac{6578}{147483}\right)\cdot g^{1,1}(u).$$
}}

If we choose $u_0=3$ we have
\begin{itemize}
\item it is easy to see  that if $h(u)$ is a nonconstant divisor of
rad$_{\mathcal O_K[u]}((e_1(u)-e_2(u))\cdot(e_1(u)-e_3(u))$ or
rad$_{\mathcal O_K[u]}(e_2(u)-e_1(u))\cdot(e_2(u)-e_3u)))$ or
rad$_{\mathcal O_K[u]}((e_3(u)-e_1(u))\cdot(e_3(u)-e_2(u))$ in $\mathcal O_K[u]$,
then $h(3)$ is not a square in $ K.$ Thus Theorem \ref{mainh} is
satisfied for $K=\mathbb Q(q)$, so we conclude that the
specialization homomorphism $\sigma_{3}^K:E^{1,1}_g(K(u))\rightarrow
E^{1,1}_g(3)(K)$ is an injection. After taking its restriction to
$E^{1,1}_g(\mathbb Q(u))$ we conclude that
$\sigma_{3}=\sigma_{3}^K\vert_{E^{1,1}_g(\mathbb
Q(u))}:E^{1,1}_g(\mathbb Q(u))\rightarrow E^{1,1}_g(3)(\mathbb Q)$
is also an injection.
\item Magma showed that
$E^{1,1}_g(3)(\mathbb Q)$ has rank $2$.
\end{itemize}

Thus we have shown by using Lemma \ref{submain} and Lemma \ref{GL}
that
 $P^{1,1}_g,Q^{1,1}_g$ are the generators of $E^{1,1}_g(\mathbb Q(u))$ which has rank 2.

\qed

\begin{remark}

\begin{itemize}

\item[(i)] For obtaining the above table for some values of $a$ and $b$ (and
moreover proving the claim of Theorem \ref{corall} for a wider
class of values $a,b$) we observed integer values $u_0=1,\ldots 70$ (minimal) that satisfy Theorem \ref{mainh} or Theorem
\ref{tmtwist} respectively,  as well as that the root number of $
E^{a,b}_g(u_0)$ is 1 and after that  we let {\tt{ mwrank}} or
Magma try to calculate the rank. For the last two columns we stopped the search for specific values $a,b$ after we came to a good candidate $u_0$ for which we couldn't calculate the rank, while in the first two we continued the search despite. We mention that in the rational case only for $(r_1,r_2,a,b)=(4,9,-133,468)$ we took a rational $u_0=\frac {25}7$  and for the quadratic case  for $(q,r_1,r_2,a,b)=(6,2,3,-66,-200)$ we took $u_0=\frac {9}8$ and for $(q,r_1,r_2,a,b)=(17,3,1,-44,-48)$ we took $u_0=\frac 74$.
\item[(ii)] The first main column in the above table is for the rational case in a) and it lists the values for $r_1=1,2,3$ and $r_2=r_1,r_1+1,\ldots, 15$. The second main column is for the quadratic case in b) and it lists at least two examples for each quadratic field observed (for all $q\ne -3$ we take $(r_1,r_2)=(1,1),(1,2)$).
The third main column is for the cyclic cubic case in c), it was obtained using the family of curves $y^2=x^3-mx^2+(mn-3n^2)x+n^3$, where $1\leq m,n\leq 15$ relatively prime (sorted my $m$). And the last main column is for the symmetric case in d).
\item[(iii)] In the case b) of the Theorem we omit the case $r_1=r_2$ since then $a=0$ and so $g=0$, and we omit $(r_1,r_2)=(1,3)$ since then for $u_0=1,\ldots, 70$  we were not able with {\tt mwrank } to calculate the rank when the Theorem was satisfied and the root number was one.
\item[(iv)] If an element is a square in a number field  was checked using the command {\tt nfroots} and the Galois group using the command  {\tt polgalois} in Pari.
We also adjusted the function {\tt factornf} to factor in $\mathcal O_K[u]$ instead of $K[u]$.
\end{itemize}
\end{remark}

The following is the case when the splitting field of $E^{a,b}$ is  of class number two, specifically $K=\mathbb {Q}(\sqrt{-5})$, so we apply Theorem \ref{mainh} to find a peculiar injective specialization $\sigma_{u_0}$. Thus we have to know something about $\mathcal R_K=\mathcal R_{\mathbb Q(\sqrt{-5})}$.  By the example in \cite[p.129]{Knapp} we know that $\mathcal R_{\mathbb Q(\sqrt{-5})}$ is the localization of $\mathcal O_K=\mathcal O_{\mathbb Q(\sqrt{-5})}$ by the multiplicative set $S=\{1,2,2^2,2^3,2^4,2^5,\ldots\}$, where the group of units is generated by $\mathcal O_{\mathbb Q(\sqrt{-5})}^{\times}=\{\pm 1\}$ and $2$.\\
So $K=\mathbb Q(\sqrt{-5})$, $O_K=\mathbb Z[\sqrt{-5}]$ and $\mathcal R_K=S^{-1}\mathcal O_K$. For each ideal $I$ of $\mathcal O_K$ let us define $I_S:=S^{-1}I$. Then $I_S$ is an ideal of $\mathcal R_K$ and $I_S$ is proper if and only if $I\cap S=\emptyset$. The non-zero prime ideals of $\mathcal R_K$ are exactly $S^{-1}I$ where
$I$ goes through non-zero prime ideals of $\mathcal O_K$ different from $\mathcal{P}=(2,1-\sqrt{-5})$ (note that $\mathcal{P}$ is not principal in $\mathcal O_K$ and that  $(2)=\mathcal{P}^2$ in $\mathcal O_K$). Namely, if $\mathcal{I}$ is a nonzero prime ideal in $\mathcal R_K$ then $\mathcal{I}\cap \mathcal O_K$ is a nonzero prime ideal of $\mathcal O_K$ and $\mathcal{I}=S^{-1}(\mathcal{I}\cap \mathcal O_K)$. We see that
$\mathcal{I}\cap \mathcal O_K\neq \mathcal{P}$, contrary $1\in \mathcal{I}$. On the other hand, let $\mathcal{Q}$ be any nonzero prime ideal of $\mathcal O_K$ different from $\mathcal{P}$. Then $S^{-1}\mathcal{Q}$ is a prime ideal of $\mathcal R_K$. We have only to see that $S^{-1}\mathcal{Q}\neq \mathcal R_K$. Contrary, $1\in S^{-1}\mathcal{Q}$, hence $2\in \mathcal{Q}$, which implies $\mathcal{Q}|(2)$ in $\mathcal O_K$, so $\mathcal{Q}=\mathcal{P}$.\\
In the following theorem we need the decomposition of the ideal $(3)_S$ into the product of two principal ideals in $\mathcal R_K$. Recall that $(3)=\mathcal{Q}\bar{\mathcal{Q}}$ in $\mathcal O_K$, where $\mathcal{Q}=(3,1+\sqrt{-5}),\ \bar{\mathcal{Q}}=
(3,1-\sqrt{-5})$ are nonprincipal in $\mathcal O_K$. We claim that $\mathcal{P}\bar{\mathcal{Q}}=(1-\sqrt{-5})$. Namely,
$\mathcal{P}\bar{\mathcal{Q}}=(2,1+\sqrt{-5})(3,1-\sqrt{-5})=(6,2-2\sqrt{-5},3+3\sqrt{5})
=(6,2-2\sqrt{-5},3-3\sqrt{5})= (1-\sqrt{-5})$. Now we have $(1-\sqrt{-5})_S=\mathcal{P}_S\bar{\mathcal{Q}}_S=
\bar{\mathcal{Q}}_S$, especially $1\pm \sqrt{-5}$ are irreducible in $\mathcal R_K$. Note, also, that
$2-\sqrt{-5}=-1/2(1+\sqrt{-5})^2$ in $\mathcal R_K$.

Similarly one shows that $3\pm \sqrt{-5}$ and $9\pm \sqrt{-5}$ are irreducible in $\mathcal R_K$.
We have the following factorization into prime ideals in $\mathcal O_K$: $(3-\sqrt{-5})=(7,3-\sqrt{-5})(2,1-\sqrt{-5})$ and  $(9-\sqrt{-5})=(2,1-\sqrt{-5})(43,9-\sqrt{-5})$. Now after localization we get $(3-\sqrt{-5})_S=(7,3-\sqrt{-5})_S(2,1-\sqrt{-5})_S=(7,3-\sqrt{-5})_S$ and $(9-\sqrt{-5})_S=(43,9-\sqrt{-5})_S$ in $\mathcal R_K$. So they are irreducible. Also $\sqrt{-5}$, $6\pm \sqrt{-5}$ and $11$ are irreducible, since $(\sqrt{-5})$, $(6\pm \sqrt{-5})$ and $(11)$ are prime ideals in $\mathcal O_K$ and are disjoint with $S$.

Since $K$ is the quotient field of the unique factorization domain $\mathcal R_K$, thus we can obtain the irreducible nonconstant factors of a polynomial in $\mathcal R_K[u]$ by observing the factorization in $K[u]$.

\begin{theorem} \label{corh2}
Let $(a,b)=(2,12), (-7,36)$ or $(-22,84)$.
Let $E^{a,b}$ be the elliptic curve given with the equation over
$\mathbb Q$
$$E^{a,b}:y^2=x^3+ax+b.$$
Then the elliptic curve $E^{a,b}_g$   has rank
two over $\mathbb Q(u)$, with free generators the two points $P^{a,b}_g$ and
$Q^{a,b}_g$.
\end{theorem}

\pf

Here $K=\mathbb Q (\sqrt{-5})$.

{\bf First we explain in detail the case $(a,b)=(2,12)$:}

\medskip

We have
$$E^{2,12}:y^2=(x-(-2))(x-(1+\sqrt{-5}))(x-(1-\sqrt{-5})),$$
and
$g^{2,12}(u)=-2^6\cdot 3\cdot  (u^2+1)(3u^4+2u^2+2)(3u^4+4u^2+3)(2u^4+2u^2+3).$
Thus  we look at the elliptic curve $E^{a,b}_g$ over $K(u)$ given by the equation
$$E^{2,12}_g:y^2=(x-(-2)\cdot g^{2,12}(u))(x-(1+\sqrt{-5})\cdot g^{2,12}(u))(x-(1-\sqrt{-5})\cdot g^{2,12}(u)),$$
where $K=\mathbb Q(\sqrt{-5})$ is of class number two and further
{\small{
$$e_1(u)=-2\cdot g^{2,12}(u),\ \ \ e_2(u)=(1+\sqrt{-5})\cdot g^{2,12}(u),\
\ \ e_3(u)=(1-\sqrt{-5})\cdot g^{2,12}(u),$$ and so
}}
\noindent
 {\small{
$(e_1(u)-e_2(u))\cdot(e_1(u)-e_3(u))\cdot (e_2(u)-e_3(u))=$
$$=(3+\sqrt{-5})\cdot (3-\sqrt{-5})\cdot (2\cdot \sqrt{-5})\cdot (g^{2,12}(u))^3=$$
{\tiny{
$$=-\sqrt{-5}\cdot (3+\sqrt{-5})(3-\sqrt{-5})\cdot 2^{19}\cdot 3^3\cdot (u^2+1)^3(3u^4+2u^2+2)^3(3u^4+4u^2+3)^3(2u^4+2u^2+3)^3,$$
}}
which we have to factor into irreducibles in $\mathcal R_K[u]$, only the radical is of importance to get the square-free factors in $\mathcal R_K[u]$.
Above we showed that  $\sqrt{-5},1\pm \sqrt{-5}, 3\pm \sqrt{-5}$ are irreducible elements  in the principal ideal domain $\mathcal{R}_{K}=\mathcal{R}_{\mathbb Q(\sqrt{-5})}$. $2$ and $-1$ are invertible elements in $\mathcal R_K$. We also have $3=\frac 12(1+\sqrt{-5})(1-\sqrt{-5}),\ 2\pm \sqrt{-5}=-\frac 12(1\mp\sqrt{-5})^2.$

We first factor in $K[u]$:
$$(u^2+1)(3u^4+2u^2+2)(3u^4+4u^2+3)(2u^4+2u^2+3)=$$
{\tiny{
$$=2\cdot 3^2\cdot (u^2 + 1)   \left(u^2 + \frac{1+\sqrt{-5}}2\right)  \left(u^2 +\frac {  1-\sqrt{-5}}2\right) \left(u^2 + \frac {1+\sqrt{-5} }3\right)   \left(u^2 + \frac{  1-\sqrt{-5}}3\right) \left(u^2 + \frac{2+\sqrt{-5}}3\right) \left(u^2 + \frac{2-\sqrt{-5} }3\right)    $$
$$=2 \cdot(u^2 + 1)\cdot \left(u^2 + \frac{1+\sqrt{-5}}2\right)  \left(u^2 +\frac {  1-\sqrt{-5}}2\right)    \left(\frac  {1-\sqrt{-5} }2u^2 +  1\right)  \left(\frac  {1+\sqrt{-5} }2u^2 +1\right)  \left((1+\sqrt{-5})u^2 - (1-\sqrt{-5})\right)   \left((1-\sqrt{-5} )u^2 -(1+\sqrt{-5})\right)$$
}}

Now it is easy to see that the radical in $\mathcal R_K[u]$ of $(e_1(u)-e_2(u))\cdot(e_1(u)-e_3(u))\cdot (e_2(u)-e_3(u))$ factors into irreducible elements in the UFD $\mathcal R_K[u]$ as
$$\mbox{rad}_{\mathcal R_K[u]}[(e_1(u)-e_2(u))\cdot(e_1(u)-e_3(u))\cdot (e_2(u)-e_3(u))]=$$
{\tiny{
$$= \sqrt{-5}\cdot (1+\sqrt{-5})\cdot (1-\sqrt{-5})\cdot (3+\sqrt{-5})\cdot (3-\sqrt{-5})\cdot (u^2 + 1)\cdot$$
$$\cdot   \left(u^2 + \frac{1+\sqrt{-5}}2\right)  \left(u^2 +\frac {  1-\sqrt{-5}}2\right)    \left(\frac  {1-\sqrt{-5} }2u^2 +  1\right)  \left(\frac  {1+\sqrt{-5} }2u^2 +1\right) \cdot$$
$$\cdot \left((1+\sqrt{-5})u^2 - (1-\sqrt{-5})\right)   \left((1-\sqrt{-5} )u^2 -(1+\sqrt{-5})\right) .$$
}}
So we obtain all nonconstant square-free divisors of $(e_1(u)-e_2(u))\cdot(e_1(u)-e_3(u))\cdot (e_2(u)-e_3(u))$ in $\mathcal R_K[u]$ as nonconstant elements
{\tiny{
$$(-1)^{i_1}\cdot 2^{i_2}\cdot \sqrt{-5}^{i_3}(1+\sqrt{-5})^{i_4}\cdot (1-\sqrt{-5})^{i_5}\cdot (3+\sqrt{-5})^{i_6}\cdot (3-\sqrt{-5})^{i_7}\cdot  (u^2 + 1)^{i_8}\cdot  \left(u^2 + \frac{1+\sqrt{-5}}2\right)^{i_9} \cdot  \left(u^2 +\frac {  1-\sqrt{-5}}2\right)^{i_{10}}\cdot$$
$$\cdot    \left(\frac  {1-\sqrt{-5} }2u^2 +1\right)^{i_{11}} \left(\frac  {1+\sqrt{-5} }2u^2 + 1\right)^{i_{12}}    \left((1+\sqrt{-5})u^2 - (1-\sqrt{-5})\right)^{i_{13}}   \left((1-\sqrt{-5} )u^2 -(1+\sqrt{-5})\right)^{i_{14}}  ,$$
}}\noindent
where $i_k\in\{0,1\}$, for $k=1,\ldots 14.$  We had to take into account  $-1$ and $2$, the two generators of the group of units in $\mathcal R_K$.

If we choose $u_0=4$ we have
\begin{itemize}
\item  it is easy to see  that if $h(u)$ is a nonconstant square-free divisor of rad$_{\mathcal R_K[u]}((e_1(u)-e_2(u))\cdot(e_1(u)-e_3(u)))\cdot(e_2(u)-e_3(u)))$ in $\mathcal R_K[u]$, then $h(4)$ is not a square in $K.$
Thus Theorem \ref{mainh} is satisfied for $K=\mathbb Q(\sqrt{-5})$ which has class number two, so we conclude that the
specialization homomorphism $\sigma_{u_0}^K:E^{2,12}_g(K(u))\rightarrow E^{2,12}_g(u_0)(K)$ is an injection. So certainly also the specialization $\sigma_{u_0}=\sigma_{u_0}^K|_{E^{2,12}_g(\mathbb Q(u))}:E^{2,12}_g(\mathbb Q(u))\rightarrow E^{2,12}_g(u_0)(\mathbb Q)$ is an injection.

\item {\tt {MordellWeilShaInformation}}  showed that  $E^{2,12}_g(4)$ has rank $2$ over $\mathbb Q$.
\end{itemize}
Thus for   $(a,b)=(2,12)$ and  the chosen specialization homomorphism $\sigma_{u_0}$, Lemma \ref{submain} is satisfied, so we conclude
that the rank of $E^{2,12}_g$ over $\mathbb Q(u)$ is two. Now by
applying Lemma \ref{GL}  we conclude that $E^{2,12}_g$ over
$\mathbb Q(u)$ has rank two and free generators
$P^{2,12}_g,Q^{2,12}_g$.

\medskip

{\bf In short the case $(a,b)=(-7,36)$:}

We have
$$E^{-7,36}:y^2=(x-(-4))(x-(2+\sqrt{-5})(x-(2-\sqrt{-5}))=x^3-7x+36.$$

Since
$$g^{-7,36}(u)=2^2 \cdot 3^2 \cdot 7  \cdot (u^2 + 1)  (3u^2 - 4u + 3)  (3u^2 + 4u + 3)  (9u^4 + 16u^2 + 16) (16u^4 + 16u^2 + 9),$$
we have
\noindent
$(e_1(u)-e_2(u))\cdot(e_1(u)-e_3(u))\cdot (e_2(u)-e_3(u))=$
$$=(6+\sqrt{-5})\cdot (6-\sqrt{-5})\cdot (2\cdot \sqrt{-5})\cdot (g^{-7,36}(u))^3=$$
 {\tiny{
$$=\sqrt{-5}\cdot (6+\sqrt{-5})(6-\sqrt{-5})\cdot 2^7\cdot 3^6\cdot 7^3\cdot (u^2 + 1)^3  (3u^2 - 4u + 3)^3  (3u^2 + 4u + 3)^3  (9u^4 + 16u^2 + 16)^3 (16u^4 + 16u^2 + 9)^3
.$$}}
Since $\sqrt {-5}$, $1\pm \sqrt{-5}$, $3\pm\sqrt{-5}$, $6\pm\sqrt{-5}$  are irreducible in the principal ideal domain $\mathcal{R}_{K}=\mathcal{R}_{\mathbb Q(\sqrt{-5})}$  and we have
{\scriptsize{
$3=\frac 12(1-\sqrt{-5})(1+\sqrt{-5}),\ 7=\frac 12(3-\sqrt{-5})(3+\sqrt{-5}), 2\pm\sqrt{-5}=-\frac 12(1\mp\sqrt{-5})^2,$
}}then
$$\mbox{rad}_{\mathcal R_K[u]}[(e_1(u)-e_2(u))\cdot(e_1(u)-e_3(u))\cdot (e_2(u)-e_3(u))]=$$
{\tiny{
$$= \sqrt{-5}\cdot (6+\sqrt{-5})(6-\sqrt{-5})(1-\sqrt{-5})(1+\sqrt{-5})(3-\sqrt{-5})(3+\sqrt{-5})\cdot  (u^2 + 1)\cdot$$
$$\cdot  \left((1+\sqrt{-5})u + (-1+\sqrt{-5})\right)  \left((1-\sqrt{-5})u +(-1-\sqrt{-5})\right) \left((1-\sqrt{-5})u + (1+\sqrt{-5})\right)\left((1+\sqrt{-5})u + (1-\sqrt{-5}\right)   \cdot$$
$$\cdot  \left(u^2 -\frac 18(1+\sqrt{-5})^2\right) \left(u^2 -\frac 18(1-\sqrt{-5})^2\right) \left((1+\sqrt{-5})^2u^2 -8\right) \left((1-\sqrt{-5})^2u^2 -8\right)    .$$
}}

The adequate specialization is for $u_0=16$.

{\bf In short the case $(a,b)=(-22,84)$:}

We have
{\small{
$$E^{-22,84}:y^2=(x-(-6))(x-(3+\sqrt{-5})(x-(3-\sqrt{-5}))=x^3-22x+84.$$
}}
So
$$g^{-22,84}(u)=2^6 \cdot 3 \cdot7 \cdot 11 \cdot (u^2 + 1)  (7u^4 - 4u^2 + 7) (7u^4 + 18u^2 + 18)  (18u^4 + 18u^2 + 7),$$
and further
\noindent
 {\small{
$(e_1(u)-e_2(u))\cdot(e_1(u)-e_3(u))\cdot (e_2(u)-e_3(u))=$
$$=(9+\sqrt{-5})\cdot (9-\sqrt{-5})\cdot (2\cdot \sqrt{-5})\cdot (g^{-22,84}(u))^3=$$
{\tiny{
$$=\sqrt{-5}\cdot (9+\sqrt{-5})(9-\sqrt{-5})\cdot 2^{19}\cdot 3^3\cdot 7^3\cdot 11^3\cdot (u^2 + 1)^3 (7u^4 - 4u^2 + 7)^3  (7u^4 + 18u^2 + 18)^3  (18u^4 + 18u^2 + 7)^3
.$$}}
 It can be shown $\sqrt{-5},11,1\pm \sqrt{-5}, 3\pm \sqrt{-5}, 9\pm \sqrt{-5}$ are irreducible in the principal ideal domain $\mathcal{R}_{K}=\mathcal{R}_{\mathbb Q(\sqrt{-5})}$ (see above) and
{\scriptsize{
 $3=\frac 12(1-\sqrt{-5})(1+\sqrt{-5}),\ 7=\frac 12(3-\sqrt{-5})(3+\sqrt{-5}), \ 2\pm 3\sqrt{-5}=\frac 12(3\pm\sqrt{-5})^2.$}}

Now it is easy to see that the radical in $\mathcal R_K[u]$ is

$$\mbox{rad}_{\mathcal R_K[u]}[(e_1(u)-e_2(u))\cdot(e_1(u)-e_3(u))\cdot (e_2(u)-e_3(u))]=$$
{\tiny{
$$= \sqrt{-5}\cdot (9+\sqrt{-5})(9-\sqrt{-5})(1-\sqrt{-5})(1+\sqrt{-5})(3-\sqrt{-5})(3+\sqrt{-5})\cdot 11\cdot  $$
$$\cdot (u^2 + 1)  ((3+\sqrt{-5})u^2 + (-3+\sqrt{-5}))((3-\sqrt{-5})u^2 + (-3-\sqrt{-5})) \cdot$$
$$\cdot ((3-\sqrt{-5})u^2 +(1-\sqrt{-5})(1+\sqrt{-5}))((3+\sqrt{-5})u^2 +(1-\sqrt{-5})(1+\sqrt{-5}))\cdot$$
$$\cdot ((1-\sqrt{-5})(1+\sqrt{-5})u^2 + (3+\sqrt{-5})) ((1-\sqrt{-5})(1+\sqrt{-5})u^2 + (3-\sqrt{-5}))     .$$
}}

The adequate specialization is for $u_0=4$.

\qed

By comparing the speed of the programs written in Pari for
checking the conditions of Theorem \ref{mainh} compared to that of
Theorem \ref{tmtwist} for the quadratic case, we
conclude that Theorem \ref{tmtwist} (i) is much faster. We have to mention that Theorem \ref{mainh} gives the result
for the specialization $\sigma_{u_0}^K:E^{a,b}_g(K(u))\rightarrow
E^{a,b}_g(u_0)(K)$ where $K$ is the corresponding splitting field of $E^{a,b}$, while Theorem \ref{tmtwist} deals with
$\sigma_{u_0}:E^{a,b}_g(\mathbb Q(u))\rightarrow
E^{a,b}_g(u_0)(\mathbb Q)$.

Theorems \ref{main}, Theorem \ref{mainh} and Theorem \ref{tmtwist} in general don't detect all
injective specializations.

These theorems were used to find an
adequate injective specialization homomorphism and with it to prove
that $P^{a,b}_g$ and $Q_g^{a,b}$ are the free generators of
$E^{a,b}_g$ over $\mathbb Q(u)$,  for each of the values $a,b$ in
the above Theorem \ref{corall} and Theorem \ref{corh2}.

So for each $a,b$ from Theorem  \ref{corall} we can now check
for some range of values $u_0$ (for example  $u_0=1,\ldots,80$)
which exactly of the corresponding specializations
$\sigma_{u_0}:E^{a,b}_g(\mathbb Q(u))\rightarrow
E^{a,b}_g(u_0)(\mathbb Q)$ are injective and which not in the
following way. We can do this because we know  that the
specialization homomorphisms are injections on the torsion part and
we know the free generators of $E^{a,b}_g$ over $\mathbb Q(u)$, so
we can apply Lemma \ref{inj} and translate the question of the
injectivity of the specialization  $\sigma_{u_0}$ to the question
of the independence of the specialized points $P^{a,b}_g(u_0)$ and
$Q^{a,b}_g(u_0)$ in $E^{a,b}_g(u_0)(\mathbb Q)$.  We checked in
this way for all $a,b$ from Theorem \ref{corall} precisely the
injectivity of all specializations $\sigma_{u_0}$ in the range
$u_0=1,\ldots,80$ with Magma's command {\tt {IsLinearlyIndependent}} and we obtained that for all except
the four quadratic cases
$(a,b)=(6,20),(-14,-20),(-56,-160),(-126,-540)$ only  $u_0=1$
gives a non-injective specialization. For the excluded ones  there
are two non-injective specializations in this range, these are for
$u_0=1,2$.

We mention that the theorems don't detect all of the
injections in this range. For example, in the rational case
$(a,b)=(-7,6)$ we don't get the answer about the injectivity of
$\sigma_{u_0}$ with  Theorem \ref{main} for
$u_0=1,2,3,4,5,6,7,9,10,11,12,
15,29,37,$ $38,40,41,45,46,49,54,56,58,60,71$. In the quadratic
case   $(a,b)=(1,10)$ (splitting field $\mathbb Q(\sqrt{-1})$) we
don't get the answer with Theorem \ref{tmtwist} (i) for
$u_0=1,2,3,5,7,8,14,$ $18,38,41,57.$ In the cyclic cubic case
$(a,b)=(-1647,1647)$ with Theorem \ref{tmtwist} (ii)  for
$u_0=1,2,3,4,5$ and in the symmetric case $(a,b)=(1,1)$ the only
specialization $\sigma_{u_0}:E^{a,b}_g(K(u))\rightarrow
E^{a,b}_g(u_0)(K)$ in  this range  for which we don't know if it
is injective using Theorem \ref{main} is for $u_0=1$, here $K$ is
the corresponding splitting field.

As a conclusion to this section we conjecture:
\begin{conjecture}

For all $a,b$ such that $a\cdot b\ne 0$,
 the elliptic curve $E^{a,b}_g$   has rank
two over $\mathbb Q(u)$, with free generators the two points $P^{a,b}_g$ and
$Q^{a,b}_g$.
\end{conjecture}

\medskip

{\footnotesize \noindent Ivica Gusi\'c \\
Faculty of Chemical Engin. and Techn. \\
University of Zagreb
\\ Maruli\'cev trg 19, 10000 Zagreb \\ Croatia \\
{\em E-mail address}: {\tt igusic@fkit.hr}}

\bigskip

{\footnotesize \noindent Petra Tadi\'c \\
Geotehnical faculty\\
University of Zagreb\\
Hallerova aleja 7, 42000 Vara\v zdin\\
Croatia\\
 {\em  E-mail adress}: {\tt petra.tadic.zg@gmail.com,\ ptadic@gfv.hr
 }

\end{document}